\documentclass[12pt]{article}
\usepackage[left=20mm, top=20mm, right=20mm, bottom=20mm]{geometry}
\usepackage{amsmath}
\usepackage[font=scriptsize]{caption}
\usepackage{multicol}
\usepackage{graphicx}
\usepackage{footmisc}
\usepackage{float}
\usepackage{longtable}
\usepackage{makeidx}
\usepackage{endnotes}
\usepackage[norefeq]{nomencl}
\makeglossary

\begin{document}

\title{\textbf{A modified fifth-order WENO scheme for hyperbolic conservation
laws}}
\author{Samala Rathan\thanks{Email: rathan.maths@gmail.com},\hspace{1mm}
G Naga Raju \thanks{Email: gnagaraju@mth.vnit.ac.in}\\
 {\small{}{}{}Department of Mathematics, Visvesvaraya National Institute
of Technology, Nagpur, India}}
\date{}
\maketitle
\begin{abstract}
This paper deals with a new fifth-order weighted essentially non-oscillatory
(WENO) scheme improving the WENO-NS and WENO-P methods which are introduced
in Ha et al. J. Comput. Phys. (2013) and Kim et al., J. Sci. Comput.
(2016) respectively. These two schemes provide the fifth-order accuracy
at the critical points where the first derivatives vanish but the
second derivatives are non-zero. In this paper, we have presented
a scheme by defining a new global-smoothness indicator which shows
an improved behavior over the solution to the WENO-NS and WENO-P schemes
and the proposed scheme attains optimal approximation order, even
at the critical points where the first and second derivatives vanish
but the third derivatives are non-zero.
\end{abstract}

\section{Introduction}

\hspace{0.6cm}The hyperbolic conservation laws may develop discontinuities
in its solution even if the initial conditions are smooth. Therefore
classical numerical methods which depend on Taylor's expansion fails,
so the spurious oscillations occur in the solution. To resolve this,
the Total Variation Diminishing (TVD) schemes are constructed by Harten
in \cite{Ha83,Ha84} based on the principle that, the total variation
to the approximation of numerical solution must be non-increasing
in time. But the TVD schemes are at most first order accurate near
smooth extrema \cite{OC84}.

In order to overcome this difficulty, Harten et al. \cite{HOEC86,HO87,HO97}
succeeded by relaxing the TVD condition and allowing the spurious
oscillations to occur in the numerical scheme, but the $O(1)$ Gibbs-like
phenomena are essentially prevented. This is the first successful
higher order spatial discretization method for the hyperbolic conservation
laws that achieves the essentially non-oscillatory (ENO) property
and known as ENO schemes. In \cite{HOEC86}, finite-volume ENO method
was studied and shown that, to have a uniform high-order accuracy
right up to the location of any discontinuity. Later, the finite-difference
ENO scheme was developed by Shu and Osher in \cite{osherSHU,Shu-osher1}.

The ENO method works based on the idea of choosing the interpolation
points over a stencil which avoids the initiation of oscillations
in the numerical solution. To do this, a smoothness indicator of the
solution is determined over each stencil and by using this, the smoothest
candidate stencil is chosen from a set of candidate stencils. As the
result, the ENO scheme avoids spurious oscillations near discontinuities
and obtains information from smooth regions only.

The weighted ENO (WENO) scheme is introduced by Liu et al. \cite{XDLiu},
in the finite-volume framework up to third-order of accuracy and later,
Jiang and Shu \cite{Jiangandshu7}, introduced this WENO scheme
in finite-difference framework which we refer here as WENO-JS scheme
by constructing the new smoothness indicators that measure the sum
of the normalized squares of the scaled $L_{2}$ norms of all derivatives
of local interpolating polynomials. Very high order schemes are constructed
in \cite{balsaraandshu8} based on the WENO-JS scheme which satisfies
monotonicity preserving property.

The main idea behind the WENO scheme is that it uses a convex combination
of all the ENO candidate sub-stencils in a non-linear manner and assigns
a weight to each sub-stencil between $0$ and $1$ based on its local
smoothness indicator. The basic strategy of assigning the weights
is that, by the combining of the lower order polynomials with optimal
weights, which yields an upwind scheme of maximum possible order in
smooth regions of the solution and assigns smaller weights to those
lower order polynomials whose stencils contain discontinuities so
that the essentially non-oscillatory property is achieved. The detailed
description of the methodology about ENO and WENO schemes and their
implementation can be found in \cite{shunotes,cwshu16}.

It is first pointed out by Henrick et al. in \cite{henrickaslampowers5},
that the desired convergence rate of the fifth-order WENO-JS was not
achieved for many problems where the first and third derivatives of
the flux do not simultaneously vanish and the behavior of the scheme
is sensitive to $\epsilon,$ a parameter is introduced in order to
avoid division by zero in the evaluation of non-linear weights. They
suggested an improved version of the WENO-JS scheme which is termed
as mapped WENO and abbreviated by WENO-M. By using a mapping function
on the nonlinear weights, WENO-M satisfies the sufficient condition
where WENO-JS fails and obtains an optimal order of convergence near
simple smooth extrema. Consequently, higher order WENO schemes were
developed based on mapping function in \cite{Gerolymus}.

In a different approach to the construction of the nonlinear weights,
Borges et al. \cite{Borgescaramona10}, introduced the fifth-order
WENO-Z scheme. In their study, the authors measured the smoothness
of the large stencil which comprises all sub-stencils and incorporated
this in devising the smoothness indicators and nonlinear weights.
The resulting WENO-Z scheme is less dissipative than WENO-JS and WENO-M.
The convergence order of the WENO-Z scheme is four at the first-order
critical points and degrade to two when higher order critical points
are encountered. Further, a closed-form formula was derived in \cite{Castroetal}
for the WENO-Z scheme to the all odd orders, higher than fifth-order
accuracy.

The smoothness indicators based on Lagrange interpolation polynomials
are derived in \cite{Pfanetal14} which gives the desired order
of convergence if the first and second derivatives vanish but the
third derivatives are non-zero by constructing the higher-order global
smoothness indicator in $L_{2}-$sense. The resulting scheme shows
less dissipative nature than other schemes and subsequently high order
schemes were presented in \cite{Pfan15}. Modified smoothness indicators
of WENO-JS scheme is presented in \cite{MWENO-Z} based on the linear
combination of second order derivatives over the global stencil by
using the idea of WENO-Z scheme. Acker et al. in \cite{Acker}, have
observed that the improving the weights information where the solution
is non-smooth is more important than the improving the accuracy at
critical points. To do this, they introduced an extra term to the
WENO-Z weights so that the scheme behaves with the same stability
and sharpness as the WENO-Z scheme at discontinuities and shocks with
a higher numerical resolution. The WENO methodology is still in development
to improve its rate of convergence in smooth regions and decrease
the dissipation near the discontinuities even it is successful in
a wide number of applications.

Recently Ha et al. \cite{Haetal13} introduced the local as well
as global smoothness indicators based on $L_{1}$-norm approach and
termed this scheme as WENO-NS. The WENO-NS scheme provides an improved
behavior compared to other fifth-order WENO schemes for the problems
which contain discontinuities. It is well known that the smoothness
indicators constructed based on $L_{1}-$norm approach may lead to
provide loss of regularity of the solution. To overcome this difficulty,
the authors constructed local smoothness indicators by developing
an approximation method to derivatives with higher accuracy and introduce
a parameter $\xi$ in calculating the local smoothness indicators
to the tradeoff between the accuracy around the smooth and discontinuous
regions. The global smoothness indicator is constructed as an average
of the smoothness of a global stencil data and a middle stencil information
through a mapping function which satisfies desirable sufficient condition
so that the scheme achieves required order of accuracy.

The main difficulty in WENO-NS scheme is to find global smoothness
indicator and because of the symmetry nature of the two local smoothness
indicators in the WENO-NS scheme out of three, the given three sub-stencils
may provide the unbalanced contribution to the evaluation of flux
at the interface. Kim et al. \cite{KIMetal}, made a balanced tradeoff
among all the sub-stencils by introducing a parameter $\delta$ in
the formulation of local smoothness indicators and constructed a global
smoothness indicator which does not include an extra information like
as in WENO-NS scheme i.e., middle stencil data. This modification
yields better performance than the WENO-NS scheme and termed as WENO-P
scheme.

A simple analysis verifies that the resulting WENO-NS and WENO-P schemes
have the fifth-order accuracy if the first-order critical point vanish
but the second derivatives are non-zero. The main aim of this study
is to further improve the WENO-NS and WENO-P schemes to achieve desired
order of accuracy even if the first and second derivatives vanish
but the third derivative is non-zero.

In this article, we analyzed the fifth-order WENO scheme with the
smoothness indicators developed in \cite{Jiangandshu7,henrickaslampowers5,Borgescaramona10,Haetal13,KIMetal}
and derived a new global smoothness indicator which is a linear combination
of second order derivatives leads to give a fourth-order of accuracy
so, the resulting global smoothness indicator provides a much smoother
information to the evaluation point in evaluating the non-linear weights.
It is verified that the proposed WENO scheme has the fifth-order accuracy
even at critical points where first and second derivatives vanish
but the third derivative is non-zero through the Taylor expansion.
We call this scheme as modified WENO-P (MWENO-P) which takes almost
the same computational cost as that of WENO-NS or WENO-P and is simple
to implement as the WENO-JS or WENO-Z schemes. The numerical experiments
are shown that the proposed scheme MWENO-P performs better than WENO-JS,
WENO-Z, WENO-NS and WENO-P for the problems which contain discontinuities.

The organization of the paper is as follows. Preliminaries to understand
about WENO reconstructions to the one-dimensional scalar conservation
laws are presented in section 2 and in section $3$ details about
the construction of a new global measurement which estimates smoothness
of a local solution in the construction of a fifth-order WENO scheme.
Some numerical results are provided in Section $4$ to demonstrate
advantages of the proposed WENO scheme. Finally, concluding remarks
are given in Section $5$.

\section{WENO schemes}

\hspace{1cm}Consider the general form of conservation law
\begin{equation}
u_{t}+f(u)_{x}=0,-\infty<x<\infty,t>0,\label{eq:ex1}
\end{equation}
with initial condition
\[
u(x,0)=u_{0}(x).
\]
Here the function $u=(u_{1},u_{2},......,u_{m})^{T}$ is a $m$-dimensional
vector of conserved variables, flux $f(u)$ is a vector-valued function
of $m$ components, $x$ and $t$ denote space and time variables
respectively. The system is called hyperbolic if all the eigen values
$\lambda_{1},\lambda_{2},...,\lambda_{m}$ of the Jacobian matrix
$A=\frac{\partial f}{\partial u}$, of the flux function are real
and set of right eigen vectors are complete.

Let $\left\{ I_{j}\right\} $ be a partition of a spatial domain with
the $j^{th}$ cell $I_{j}=[x_{j-\frac{1}{2}},x_{j+\frac{1}{2}}]$,
where $x_{j\pm\frac{1}{2}}$ are called cell interfaces. The centre
of each cell $I_{j}$ is denoted by $x_{j}=\frac{1}{2}(x_{j+\frac{1}{2}}+x_{j-\frac{1}{2}})$
and the value of the function at the node $x_{j}$ is denoted by $f_{j}=f(x_{j})$.
We assume that the set ${\{x_{j+\frac{1}{2}}}\}_{j}$ is uniformly
spatial gridded throughout the domain and the cell size of $I_{j}$
is denoted by $\triangle x=x_{j+\frac{1}{2}}-x_{j-\frac{1}{2}}$.
The approximation of one-dimensional hyperbolic conservation laws
\eqref{eq:ex1} leads to system of ordinary differential equations
by applying the method of lines, where the finite difference approximation
is replaced to the spatial derivative and yields a semi discrete scheme
\begin{equation}
\frac{d{u_{j}}}{dt}=-\frac{1}{\triangle x}({\hat{f}_{j+\frac{1}{2}}-\hat{f}_{j-\frac{1}{2}}}).\label{2}
\end{equation}

Here ${u_{j}}$ is the approximation to the point value $u(x_{j},t)$
and $\hat{f}_{j\pm\frac{1}{2}}$ are called numerical fluxes which
are Lipschitz continuous in each of its arguments and is consistent
with the physical flux $\hat{f}(u,....,u)=f(u)$. The conservation
property is obtained by defining a function $h(x)$ implicitly through
the following equation (see Lemma $2.1$ of \cite{Shu-osher1})
\begin{equation}
f(u(x,.))=\frac{1}{\triangle x}\intop_{x_{j-\frac{1}{2}}}^{x_{j+\frac{1}{2}}}h(\xi)d\xi.\label{3}
\end{equation}
Differentiating (\ref{3}) with respect to $x$ yields
\begin{equation}
f(u(x,.))_{x}=\frac{1}{\Delta x}(h(x_{j+\frac{1}{2}})-h(x_{j-\frac{1}{2}})),\label{4}
\end{equation}
where $h(x_{j\pm\frac{1}{2}})$ is a approximation to the numerical
flux $\hat{f}_{j\pm\frac{1}{2}}$ with a high order of accuracy, that
is,
\begin{equation}
\hat{f}_{j\pm\frac{1}{2}}=h(x_{j\pm\frac{1}{2}})+O(\Delta x^{r}).\label{5}
\end{equation}
To ensure the numerical stability and to avoid entropy violating solutions,
the flux $f(u)$ is splitted into two parts $f^{+}$ and $f^{-}$,
thus
\begin{equation}
f(u)=f^{+}(u)+f^{-}(u),\label{6}
\end{equation}
where $\frac{df^{+}(u)}{du}\geq0$ and $\frac{df^{-}(u)}{du}\leq0$.

The numerical fluxes $\hat{f}_{j+\frac{1}{2}}^{+}$ and $\hat{f}_{j+\frac{1}{2}}^{-}$
evaluates at $x_{j+\frac{1}{2}}$ obtained from (\ref{3}) which are
positive and negative parts of $f(u)$ respectively and with this
we have $\hat{f}_{j+\frac{1}{2}}=\hat{f}_{j+\frac{1}{2}}^{+}+\hat{f}_{j+\frac{1}{2}}^{-}$.
We will only describe how $\hat{f}_{j+\frac{1}{2}}^{+}$ is approximated
because the negative part of the split flux, is symmetric to the positive
part with respect to $x_{j+\frac{1}{2}}$. For brevity, we drop the
'$+$' sign in the superscript from here onwards.

\subsection{Fifth order WENO schemes}

\hspace{1cm}The WENO schemes are designed for the approximation to
the spatial derivative to solve the hyperbolic conservation laws,
which are used to reconstruct the unknown values of a given flux function
$f$ from its known values in an essentially non-oscillatory manner.
To construct $\hat{f}_{j+\frac{1}{2}}$, the classical fifth-order
WENO scheme uses five point stencil $S^{5}=\{x_{j-2},x_{j-1},x_{j},x_{j+1},x_{j+2}\}$
which is subdivided into three candidate sub-stencils $S_{0}(j)$,\,$S_{1}(j)$
and $S_{2}(j)$. To each cell $I_{j}$, the corresponding stencil
is denoted by
\[
S_{k}(j)=\{x_{j+k-2},x_{j+k-1},x_{j+k}\},\,k=0,1,2
\]
and let
\[
\hat{f}_{j+\frac{1}{2}}^{k}=\sum_{q=0}^{2}c_{k,q}f_{j+k+q-2},
\]
be the second-degree polynomial constructed on the stencil $S_{k}(j)$
to approximate the value $h(x_{j+\frac{1}{2}})$ where the coefficients
$c_{k,q}(q=0,1,2)$ are the Lagrange's interpolation coefficients
depending on the shifting parameter $k$. The flux values on each
stencil can be written in the form of
\begin{eqnarray}
\hat{f}_{j+\frac{1}{2}}^{0} & = & \frac{1}{6}(2f_{j-2}-7f_{j-1}+11f_{j}),\nonumber \\
\hat{f}_{j+\frac{1}{2}}^{1} & = & \frac{1}{6}(-f_{j-1}+5f_{j}+2f_{j+1}),\label{eq:9}\\
\hat{f}_{j+\frac{1}{2}}^{2} & = & \frac{1}{6}(2f_{j}+5f_{j+1}-f_{j+2}).\nonumber
\end{eqnarray}
The flux $\hat{f}_{j-\frac{1}{2}}^{k}$ is obtained through shifting
the index by $-1$. The Taylor expansion of (\ref{eq:9}) reveals
\begin{eqnarray*}
\hat{f}_{j\pm\frac{1}{2}}^{0} & = & h_{j\pm\frac{1}{2}}-\frac{1}{4}f{}^{'''}(0)\triangle x^{3}+O(\triangle x^{4}),\\
\hat{f}_{j\pm\frac{1}{2}}^{1} & = & h_{j\pm\frac{1}{2}}+\frac{1}{12}f{}^{'''}(0)\triangle x^{3}+O(\triangle x^{4}),\\
\hat{f}_{j\pm\frac{1}{2}}^{2} & = & h_{j\pm\frac{1}{2}}-\frac{1}{12}f{}^{'''}(0)\triangle x^{3}+O(\triangle x^{4}).
\end{eqnarray*}
The convex combination of these flux functions define the approximation
to the value of $h({x}_{j+\frac{1}{2}})$ which is
\begin{equation}
\hat{f}_{j+\frac{1}{2}}=\sum_{k=0}^{2}\omega_{k}\hat{f}_{j+\frac{1}{2}}^{k},\label{eq:11}
\end{equation}
where $\omega_{k}$ are the non-linear weights. If the function $h(x)$
is smooth in all the sub-stencils $S_{k}(j)$, $k=0,1,2$, we calculate
the constants $d_{k}$ such that its linear combination with ${\hat{f}}_{j+\frac{1}{2}}^{k}$
gives the fifth order convergence to $h({x}_{j+\frac{1}{2}})$, that
is,
\[
h_{j+\frac{1}{2}}=\sum_{k=0}^{2}d_{k}\hat{f}_{j+\frac{1}{2}}^{k}+O(\triangle x^{5}).
\]
The coefficients $d_{k}$ are known as the ideal weights because they
generate the upstream central fifth-order scheme for the five-point
stencil. The values of ideal weights are given by
\begin{equation}
d_{0}=1/10,\,\,d_{1}=6/10,\,\,d_{2}=3/10.\label{13}
\end{equation}
The non-linear weights $\omega_{k}$ are defined in (\ref{eq:11})
as
\begin{equation}
\omega_{k}=\frac{\alpha_{k}}{\sum_{q=0}^{2}\alpha_{q}},\label{14}
\end{equation}
and
\begin{equation}
\alpha_{k}=\frac{d_{k}}{(\epsilon+\beta_{k})^{2}},\label{15}
\end{equation}
where $0<\epsilon<<1$ is introduced to prevent the denominator becoming
zero and $\beta_{k}$ is a smoothness indicator of the flux $\hat{f}^{k}$
which measures the smoothness of a solution over a particular stencil.
The necessary and sufficient conditions to achieve the fifth-order
convergence for the WENO scheme are

\begin{align}
\sum_{k=0}^{2}(\omega_{k}^{\pm}-d_{k}) & =O(\triangle x^{6}),\label{eq:16}\\
\sum_{k=0}^{2}A_{k}(\omega_{k}^{+}-\omega_{k}^{-}) & =O(\triangle x^{3}),\label{eq:17}\\
\omega_{k}^{\pm}-{d}_{k} & =O(\triangle x^{2}).\label{eq:18}
\end{align}
Since (\ref{eq:16}) holds always due to the normalization, a sufficient
condition for the fifth-order convergence to the scheme is derived
in \cite{Borgescaramona10} as
\begin{equation}
\omega_{k}^{\pm}-{d}_{k}=O(\triangle x^{3}),\label{20}
\end{equation}
where the superscripts $'+'$ and $'-'$ on $\omega_{k}$ correspond
to their use in either $f_{j+\frac{1}{2}}^{k}$ and $f_{j-\frac{1}{2}}^{k}$
respectively.

\subsubsection{The WENO-JS scheme}

\hspace{1cm}The suggested smoothness indicators $\beta_{k}$ of Jiang
and Shu in \cite{Jiangandshu7} are given by
\begin{equation}
\beta_{k}=\sum_{q=1}^{2}\triangle x^{2q-1}\intop_{x-\frac{\triangle x}{2}}^{x+\frac{\triangle x}{2}}\left(\frac{d^{q}\hat{f}^{k}}{dx^{q}}\right)^{2}dx.\label{21}
\end{equation}
The explicit form of the these smoothness indicators are as follows:
\begin{eqnarray}
\beta_{0} & = & \frac{13}{12}(f_{j-2}-2f_{j-1}+f_{j})^{2}+\frac{1}{4}(f_{j-2}-4f_{j-1}+3f_{j})^{2},\nonumber \\
\beta_{1} & = & \frac{13}{12}(f_{j-1}-2f_{j}+f_{j+1})^{2}+\frac{1}{4}(f_{j+1}-f_{j-1})^{2},\label{eq:22}\\
\beta_{2} & = & \frac{13}{12}(f_{j}-2f_{j+1}+f_{j+2})^{2}+\frac{1}{4}(3f_{j}-4f_{j+1}+f_{j+2})^{2}.\nonumber
\end{eqnarray}
By the Taylor's expansion of these smoothness indicators, one obtain
\begin{eqnarray*}
\beta_{0} & = & f^{'^{2}}\triangle x^{2}+(\frac{13}{12}f^{''^{2}}-\frac{2}{3}f\,^{'}f^{'''})\triangle x^{4}+(\frac{-13}{6}f^{''}f^{'''}+\frac{1}{2}f\,^{'}f^{(4)})\triangle x^{5}+O(\triangle x^{6}),\\
\beta_{1} & = & f^{'^{2}}\triangle x^{2}+(\frac{13}{12}f^{''^{2}}+\frac{1}{3}f^{'}f^{'''})\triangle x^{4}+O(\triangle x^{6}),\\
\beta_{2} & = & f^{'^{2}}\triangle x^{2}+(\frac{13}{12}f^{''^{2}}-\frac{2}{3}f^{'}f^{'''})\triangle x^{4}+(\frac{13}{6}f^{''}f^{'''}-\frac{1}{2}f^{'}f^{(4)})\triangle x^{5}+O(\triangle x^{6}).
\end{eqnarray*}

The smoothness indicators (\ref{eq:22}) are smoother than the smoothness
indicators presented in \cite{XDLiu}, which are constructed based
on the total variation measurement of a stencil in $L_{1}-$norm.
The Taylor expansion of these smoothness indicators $\beta_{k}$ satisfies
the sufficient condition
\[
\beta_{k}=D(1+O(\Delta x^{2})),
\]
where the constant $D$ is independent of $k$ but depends on the
$\Delta x$, which implies that the WENO weights satisfy the condition
\[
\omega_{k}-{d}_{k}=O(\triangle x^{2}).
\]

Therefore, the final WENO reconstruction provides the fifth order
convergence in the smooth regions. However at the critical points
where the first derivative of $f$ vanishes, the convergence property
may not hold since $\beta_{k}=D(1+O(\Delta x))$ yields $\omega_{k}-{d}_{k}=O(\triangle x)$.
Further, the order of convergence is degraded to two if the first
derivative and second derivatives vanish but the third derivative
is non-zero .

\subsubsection{WENO-M}

\hspace{1cm}When the fifth-order WENO-JS scheme is used, the condition
(\ref{20}) may not hold at certain smooth extrema or near critical
points, yielding the third-order accuracy. To overcome this situation
Henrick et al. \cite{henrickaslampowers5} introduced a mapping
function $g_{k}(\omega)$ as
\begin{equation}
g_{k}(\omega)=\frac{\omega(d_{k}+d_{k}^{2}-3d_{k}\omega+\omega^{2})}{d_{k}^{2}+\omega(1-2d_{k})},k=0,1,2,\label{23}
\end{equation}
where ${d}_{k}$'s are the ideal weights and $\omega\in[0,1]$. This
function is a non-decreasing monotone function on $[0,1]$ which satisfies
the following properties
\begin{eqnarray}
\begin{aligned} & (i)\:0\leq g_{k}(\omega)\leq1,g_{k}(0)=0\textrm{ and }g_{k}(1)=1.\\
 & (ii)\:g_{k}(\omega)\approx0\,\textrm{ if}\:\omega\approx0,\:g_{k}(\omega)\approx1\,\text{ if}\:\omega\approx1.\\
 & (iii)\:g_{k}(d_{k})=d_{k},\,g_{k}^{\prime}(d_{k})=g_{k}^{\prime\prime}(d_{k})=0.\\
 & (iv)\:g_{k}(\omega)=d_{k}+O(\Delta x^{6}),\,\textrm{ if}\;\omega=d_{k}+O(\Delta x^{2})
\end{aligned}
 &  & \hspace{1mm}\label{24}
\end{eqnarray}
and with this mapping function, non-linear weights are defined as
\begin{equation}
\omega_{k}^{M}=\frac{\alpha_{k}^{M}}{\sum_{l=0}^{2}\alpha_{l}^{M}}\hspace{3mm}\text{and}\hspace{3mm}\alpha_{k}^{M}=g_{k}(\omega_{k}^{JS}),k=0,1,2,\label{25}
\end{equation}
where $\omega_{k}^{JS}$ are computed in (\ref{14}) with the WENO-JS
scheme by using the smoothness indicators defined in (\ref{21}).

\subsubsection{WENO-Z }

\hspace{1cm} Borges et al. \cite{Borgescaramona10} presented a scheme,
known as WENO-Z, by introducing global smoothness indicator. The idea
is to use whole five point stencil $S^{5}$ to define a new smoothness
indicator of higher order than the classical smoothness indicator
$\beta_{k}$, as
\begin{equation}
\tau_{5}=|\beta_{0}-\beta_{2}|.\label{26}
\end{equation}

By Taylor's expansion, the truncation error of $\tau_{5}$ is
\begin{equation}
\tau_{5}=\frac{13}{3}|f^{\prime\prime}f^{\prime\prime\prime}|\Delta x^{5}+O(\Delta x^{6}),\label{27}
\end{equation}
and the smoothness indicators are defined as
\begin{equation}
\beta_{k}^{z}=\frac{\beta_{k}+\epsilon}{\beta_{k}+\epsilon+\tau_{5}},\,\,k=0,1,2,\label{28}
\end{equation}
and finally with this local smoothness indicators, the WENO weights
are defined as
\begin{equation}
\omega_{k}^{z}=\frac{\alpha_{k}^{z}}{\sum_{q=0}^{2}\alpha_{q}^{z}},\,\alpha_{k}^{z}=\frac{d_{k}}{\beta_{k}^{z}}=d_{k}\bigg(1+\bigg[\frac{\tau_{5}}{\beta_{k}+\epsilon}\bigg]^{p}\bigg),\,k=0,1,2,\label{eq:zwt}
\end{equation}
where $\epsilon$ is chosen as a very small value in order to force
this parameter to play only its original role of not allowing vanishing
denominator for each weight and with the value $p=1$ in (\ref{eq:zwt}),
the scheme achieves the fourth order accuracy and with $p=2$, the
scheme achieves the fifth-order accuracy.

\subsubsection{The WENO-NS scheme}

\hspace{1cm} Ha et al., have introduced smoothness indicators based
on $L_{1}-$ norm in \cite{Haetal13} which measures the smoothness
of a solution on each $3$-point stencil $S_{k}(j),k=0,1,2$ by estimating
the approximate magnitude of derivatives as
\begin{equation}
\beta_{k}=\xi|L_{1,k}f|+|L_{2,k}f|,\label{eq:zi}
\end{equation}
where the operators $L_{n,k}f,k=0,1,2$ are the generalized undivided
differences, defined by
\begin{align}
|L_{1,k}f| & =(1-k)f_{j-2+k}+(2k-3)f_{j-1+k}+(2-k)f_{j+k},\label{eq:lnk1}\\
|L_{2,k}f| & =f_{j-2+k}-2f_{j-1+k}+f_{j+k}.\label{eq:lnk2}
\end{align}

The number $\xi$ is a parameter which is to balance the tradeoff
between the accuracy around the smooth regions and the discontinuous
regions. The second term $L_{2,k}f$ is the same as the ones of the
WENO-JS scheme, however, this scheme uses the absolute values while
the later uses the squared ones. The advantage with these operators
$L_{n,k}f$ is that the approximation of the derivative $\Delta x^{n}f^{(n)}$
at the point, $x_{j+\frac{1}{2}}$ with the high order of accuracy
can be obtained, that is,
\begin{equation}
L_{n,k}f=\frac{d^{n}f}{dx^{n}}(x_{j+\frac{1}{2}})+O(\Delta x^{3}).
\end{equation}
By using Theorem $3.2$ of \cite{Haetal13}, the Taylor's expansion
of $\beta_{k}'s$ are

\begin{eqnarray}
\beta_{0} & = & \xi\mid\Delta xf_{j+\frac{1}{2}}^{(1)}-\frac{23}{24}\Delta x^{3}f_{j+\frac{1}{2}}^{(3)}\mid+\mid\Delta x^{2}f_{j+\frac{1}{2}}^{(2)}-\frac{3}{2}\Delta x^{3}f_{j+\frac{1}{2}}^{(3)}\mid+O(\Delta x^{4}),\nonumber \\
\beta_{1} & = & \xi\mid\Delta xf_{j+\frac{1}{2}}^{(1)}+\frac{1}{24}\Delta x^{3}f_{j+\frac{1}{2}}^{(3)}\mid+\mid\Delta x^{2}f_{j+\frac{1}{2}}^{(2)}-\frac{1}{2}\Delta x^{3}f_{j+\frac{1}{2}}^{(3)}\mid+O(\Delta x^{4}),\label{eq:22-1}\\
\beta_{2} & = & \xi\mid\Delta xf_{j+\frac{1}{2}}^{(1)}+\frac{1}{24}\Delta x^{3}f_{j+\frac{1}{2}}^{(3)}\mid+\mid\Delta x^{2}f_{j+\frac{1}{2}}^{(2)}+\frac{1}{2}\Delta x^{3}f_{j+\frac{1}{2}}^{(3)}\mid+O(\Delta x^{4}).\nonumber
\end{eqnarray}
The non-linear weights are defined as
\begin{equation}
\omega_{k}^{NS}=\frac{\alpha_{k}^{NS}}{\sum_{q=0}^{2}\alpha_{q}^{NS}},\,\,\,\alpha_{k}^{NS}=d_{k}\bigg(1+\frac{\zeta}{(\beta_{k}+\epsilon)^{2}}\bigg),\,\,\,k=0,1,2,\label{35}
\end{equation}
where
\begin{eqnarray}
\zeta=\frac{1}{2}(|\beta_{0}-\beta_{2}|^{2}+g(L_{1,1}f)^{2}),\,\,\,g(x)=\frac{x^{3}}{1+x^{3}}\label{36}
\end{eqnarray}
so that the sufficient condition $\omega_{k}^{\pm}-{d}_{k}=O(\triangle x^{3})$
holds even first derivative vanish but the second derivative is non-zero.

\subsubsection{The WENO-P scheme}

\hspace{1cm} The drawback of WENO-NS scheme is that all the three
sub-stencils $S_{0}(j)$,\,$S_{1}(j)$ and $S_{2}(j)$ may provide
unbalanced contribution to the evaluation point $x_{j+\frac{1}{2}}$
because, the stencils $S_{1}(j)$ and $S_{2}(j)$ are symmetric with
respect to the point $x_{j+\frac{1}{2}}$ whereas $S_{0}(j)$ is not.
Therefore the interpolation process may be over-influenced to the
evaluation point from the left sub-stencils. For this, in \cite{KIMetal}
the authors improved the WENO-NS scheme by adjusting smoothness indicators
to the three sub-stencils by
\begin{equation}
\tilde{\beta_{0}}=\beta_{0},\,\,\tilde{\beta_{1}}=(1+\delta)\beta_{1},\,\,\tilde{\beta_{2}}=(1-\delta)\beta_{2}.\label{eqdelta}
\end{equation}

They also observed that the global smoothness measurement which uses
an additional contribution term which measures the regularity of a
solution over the stencil $S_{1}(j)$. Instead of using the measurement
proposed in WENO-NS scheme and to save computational cost, they proposed
WENO-P scheme with modification to the weights as
\begin{equation}
\omega_{k}^{P}=\frac{\alpha_{k}^{P}}{\sum_{q=0}^{2}\alpha_{q}^{P}},\,\,\alpha_{k}^{P}=d_{k}\bigg(1+\frac{\zeta}{(\tilde{\beta_{k}}+\epsilon)^{2}}\bigg),\,\,k=0,1,2,\label{38}
\end{equation}
where
\[
\zeta=(\beta_{0}-\beta_{2})^{2},
\]
and these weights satisfy the sufficient condition to get fifth order
accuracy when the first derivative vanishes but not the second derivative.

\section{A new WENO scheme}

\hspace{1cm}To achieve the desired fifth-order accuracy when the
first and second derivatives vanish but the third derivative is non-zero,
a modified WENO-P scheme called MWENO-P is proposed here. The Taylor's
expansion of the operators $L_{n,k}f$ defined in (\ref{eq:lnk1},
\ref{eq:lnk2}) by using Theorem $3.2$ of \cite{Haetal13} are
\begin{eqnarray}
L_{1,0}f & = & \Delta xf_{j+\frac{1}{2}}^{(1)}-\frac{23}{24}\Delta x^{3}f_{j+\frac{1}{2}}^{(3)}+\Delta x^{4}f_{j+\frac{1}{2}}^{(4)}+O(\Delta x^{5}),\nonumber \\
L_{1,1}f & = & \Delta xf_{j+\frac{1}{2}}^{(1)}+\frac{1}{24}\Delta x^{3}f_{j+\frac{1}{2}}^{(3)}+O(\Delta x^{5}),\nonumber \\
L_{1,2}f & = & \Delta xf_{j+\frac{1}{2}}^{(1)}+\frac{1}{24}\Delta x^{3}f_{j+\frac{1}{2}}^{(3)}+O(\Delta x^{5}),\nonumber \\
L_{2,0}f & = & \Delta x^{2}f_{j+\frac{1}{2}}^{(2)}-\frac{3}{2}\Delta x^{3}f_{j+\frac{1}{2}}^{(3)}+\frac{29}{384}\Delta x^{4}f_{j+\frac{1}{2}}^{(4)}+O(\Delta x^{5}),\\
L_{2,1}f & = & \Delta x^{2}f_{j+\frac{1}{2}}^{(2)}-\frac{1}{2}\Delta x^{3}f_{j+\frac{1}{2}}^{(3)}+\frac{5}{24}\Delta x^{4}f_{j+\frac{1}{2}}^{(4)}+O(\Delta x^{5}),\nonumber \\
L_{2,2}f & = & \Delta x^{2}f_{j+\frac{1}{2}}^{(2)}+\frac{1}{2}\Delta x^{3}f_{j+\frac{1}{2}}^{(3)}+\frac{5}{24}\Delta x^{4}f_{j+\frac{1}{2}}^{(4)}+O(\Delta x^{5}).\nonumber
\end{eqnarray}

The main idea is here to construct a high-order global smoothness
indicator which is to satisfy the sufficient condition (\ref{20})
if the first and second derivatives vanish but the third derivative
is non-zero. Also the use of higher global smoothness indicator incurs
the less dissipation near the discontinuities in the numerical scheme
\cite{Pfanetal14}. To construct such a global smooth measurement,
we define a variable $\eta$ as
\[
\eta=|L_{2,0}f+L_{2,2}f-2L_{2,1}f|^{2}.
\]
which is a linear combination of undivided differences of second-order
derivatives leads to give a fourth-order of accuracy, so this leads
to give a much smooth measurement than the global smoothness indicators
presented in \cite{Haetal13,KIMetal}. The Taylor's expansion
of $\eta$ yields
\begin{align*}
\eta & =|-\frac{51}{384}\Delta x^{4}f_{j+\frac{1}{2}}^{(4)}+O(\Delta x^{5})|^{2},\\
 & =\Delta x^{8}(A+O(\Delta x)^{2}),
\end{align*}
and we define the new non-linear weights as
\begin{equation}
\omega_{k}^{MP}=\frac{\alpha_{k}^{MP}}{\sum_{q=0}^{2}\alpha_{q}^{MP}},\,\,\alpha_{k}^{MP}=d_{k}\bigg(1+\frac{\eta}{(\tilde{\beta_{k}}+\epsilon)^{2}}\bigg),\,\,\,k=0,1,2.\label{40}
\end{equation}

\subsection{Convergence order at critical points}

\hspace{1cm}Here we discuss the convergence analysis of the proposed
scheme MWENO-P at the critical points, that is how the new weights
$\omega_{k}^{MP}$ are approaching the ideal weights $d_{k}$ in the
presence of critical points. First, consider that there is no critical
point and assume that $\epsilon=0$, then the local smoothness indicators
$\tilde{\beta_{k}},$ as defined in (\ref{eqdelta}) with weights
as in (\ref{40}) are of the form
\begin{equation}
\tilde{\beta_{k}}=\mid\Delta xf_{j+\frac{1}{2}}^{(1)}\mid+\mid\Delta x^{2}f_{j+\frac{1}{2}}^{(2)}\mid+O(\Delta x^{3})\label{eq:36}
\end{equation}
and the global smoothness indicator is of the form
\begin{equation}
\eta=O(\Delta x^{8})\label{eq:37}
\end{equation}
By substituting equations (\ref{eq:36}) and (\ref{eq:37}) in (\ref{40}),
the sufficient condition (\ref{20}) immediately holds so that the
scheme has fifth convergence order.

If $f_{j+\frac{1}{2}}^{'}=0$ and $f_{j+\frac{1}{2}}^{''}\neq0$,
the smoothness indicators are of the form
\[
\tilde{\beta_{k}}=|f_{j+\frac{1}{2}}^{''}|\Delta x^{2}(1+O(\Delta x))
\]
then there is a constant $D$ such that
\begin{align}
1+\frac{\eta}{{\tilde{\beta_{k}}^{2}}} & =1+D\Delta x^{4}+O(\Delta x^{5})\nonumber \\
 & =(1+D\Delta x^{4})\bigg(1+\frac{O(\Delta x^{5})}{1+D\Delta x^{4}}\bigg)\nonumber \\
 & =D_{\Delta x}(1+O(\Delta x^{5}))\label{eq:43}
\end{align}
where $D_{\Delta x}=(1+D\Delta x^{4})>0$. By substituting (\ref{eq:43})
in (\ref{40}), then the sufficient condition holds hence the scheme
has the fifth-order convergence at the first-order critical points.

Finally, if $f_{j+\frac{1}{2}}^{'}=0$, $f_{j+\frac{1}{2}}^{''}=0$
and $f_{j+\frac{1}{2}}^{'''}\neq0$, the smoothness indicators takes
the form
\[
{\tilde{\beta_{k}}}=|f_{j+\frac{1}{2}}^{'''}|\Delta x^{3}(1+O(\Delta x)),
\]
then there exist a constant $\widehat{D}$ such that
\begin{align}
1+\frac{\eta}{{\tilde{\beta_{k}}^{2}}} & =1+\widehat{D}\Delta x^{2}+O(\Delta x^{3})\nonumber \\
 & =(1+\widehat{D}\Delta x^{2})\bigg(1+\frac{O(\Delta x^{3})}{1+\widehat{D}\Delta x^{2}}\bigg)\nonumber \\
 & =\widehat{D}_{\Delta x}(1+O(\Delta x^{3}))\label{eq:43-1}
\end{align}
where $\widehat{D}_{\Delta x}=(1+\widehat{D}\Delta x^{2})>0$. By
substituting equation (\ref{eq:43-1}) in (\ref{40}), one can see
that the sufficient condition holds and so the newly defined smoothness
indicators attain the fifth-order convergence even in the case where
first and second derivatives vanish but the third derivatives are
non-zero.

\section{Numerical Results}

\hspace{1cm}For time evaluation in (\ref{2}) third-order TVD Runge-Kutta
scheme (TVD RK3),
\begin{equation}
\begin{aligned} & u^{(1)}=u^{n}+\triangle t\,{L}(u^{n}),\\
 & u^{(2)}=\frac{3}{4}u^{n}+\frac{1}{4}u^{(1)}+\frac{1}{4}\triangle t\,{L}(u^{(1)}),\\
 & u^{n+1}=\frac{1}{3}u^{n}+\frac{2}{3}u^{(2)}+\frac{2}{3}\triangle t\,{L}(u^{(2)}).
\end{aligned}
\label{eq:45}
\end{equation}
and fourth-order non-TVD Runge-Kutta schemes (RK-4),
\begin{equation}
\begin{aligned} & u^{(1)}=u^{n}+\frac{1}{2}\triangle t\,{L}(u^{n}),\\
 & u^{(2)}=u^{n}+\frac{1}{2}\triangle t\,{L}(u^{1}),\\
 & u^{(3)}=u^{n}+\triangle t\,{L}(u^{2}),\\
 & u^{n+1}=\frac{1}{3}(-u^{n}+u^{(1)}+2u^{(2)}+u^{(3)})+\frac{1}{6}\triangle t\,{L}(u^{(3)}).
\end{aligned}
\label{46}
\end{equation}
where $L$ is the spatial operator, are used in following numerical
examples for evaluating the approximate solution using the proposed
scheme MWENO-P.

The elaborate details of the schemes (\ref{eq:45},\ref{46}) can
be found in \cite{Gottieledandshu21}. The computed solution of
MWENO-P is compared with WENO-JS, WENO-Z, WENO-NS and WENO-P schemes.
All the numerical results are obtained on machine having an Intel(R)
core (TM) i$7-4700$MQ processor with $8$ GB of memory. In order
to ensure fairness in the comparison, all the schemes shared the same
subroutine calls and were compiled with the same compilation options.
The only differences between the implementation of the WENO schemes
were on the subroutines for computing the corresponding weights. In
order to compare with the classical scheme, WENO-JS, we took $\epsilon=10^{-6}$
and $\epsilon=10^{-40}$ for WENO-Z, WENO-NS, WENO-P and MWENO-P schemes
in all the following test cases. To prove the effectiveness of the
scheme we considered various examples. For convergence analysis we
first evaluated examples of scalar equation and later we have presented
results for Euler equations.

\subsection{Scalar test examples}

\hspace{1cm}The convergence analysis of the schemes is presented
here by considering linear advection and Burger's equations with various
initial profiles. Some of these initial profiles contain jump discontinuity
and in some cases, the solution in time leads to contact discontinuity
and shock. Here we've used the fifth-order flux version of WENO scheme
based on Lax-Friedrich's flux splitting technique and time step is
taken as $\triangle t\sim(\triangle x)^{\frac{5}{4}}$ so that the
fourth order non-TVD Runge-Kutta method in time is effectively fifth-order.
For all the examples in this subsection we've chosen $\xi=0.1$ to
evaluate the equation (\ref{eq:zi}) and $\delta=0.05$ to evaluate
the equation (\ref{eqdelta}).
\subsubsection{Example 1:}
Consider the transport equation
\begin{equation}
u_{t}+u_{x}=0,\,\,-1\leq x\leq1,\,\,t\geq0,\label{eq:lintransp}
\end{equation}
with the initial condition
\begin{equation}
u(x,0)=\text{sin}(\pi x),\label{48}
\end{equation}
and
\begin{equation}
u(x,0)=\text{sin}(\pi x-\frac{1}{\pi}\text{sin}(\pi x)),\label{49}
\end{equation}
to test the numerical convergence of the proposed scheme. Below in
Table \ref{table:1} and Table \ref{table:2} , the $L_{1}$ and $L_{\infty}$
errors are calculated up to time $t=2$ for WENO-JS, WENO-NS and WENO-P
schemes along with the proposed MWENO-P scheme for the initial condition
(\ref{48}). We observed that MWENO-P scheme is converges slowly to
the desired order of convergence for the approximate solution. The
initial condition (\ref{49}) is a special case to test the order
of convergence since its first derivative vanish but second derivative
is non-zero. Here too, the numerical order of convergence for the
MWENO-P scheme is more accurate when compared with WENO-JS, WENO-NS
and WENO-P schemes which are shown in Table \ref{table:3} and Table
\ref{table:4}.
\begin{table}[H]
$\qquad$%
\begin{tabular}{ccccc}
\hline
N  & WENO JS  & WENO-NS  & WENO-P  & MWENO-P\tabularnewline
\hline
10  & 4.8506e-02(---)  & 4.4448e-02(---)  & 5.4913e-02(---)  & 9.2794e-03(---)\tabularnewline
20  & 2.5414e-03(4.25)  & 2.1541e-03(4.36)  & 3.1845e-03(4.10)  & 3.0628e-04(4.92)\tabularnewline
40  & 8.9204e-05(4.83)  & 3.0544e-05(6.14)  & 5.0576e-05(5.97)  & 1.0249e-05(4.90)\tabularnewline
80  & 2.7766e-06(5.00)  & 3.1879e-07(6.58)  & 8.3516e-07(5.92)  & 3.1929e-07(5.00)\tabularnewline
160  & 8.6040e-08(5.01)  & 9.9380e-09(5.00)  & 1.3716e-08(5.92)  & 9.9414e-09(5.00)\tabularnewline
320  & 2.5528e-09(5.07)  & 3.1007e-10(5.00)  & 3.0641e-10(5.48)  & 3.1008e-10(5.00)\tabularnewline
\hline
\end{tabular}

\caption{\label{table:1}$L_{\infty}$ errors of linear advection Eq. (\ref{eq:lintransp})
with initial condition (\ref{48})}
\end{table}

\begin{table}[H]
$\qquad$%
\begin{tabular}{ccccc}
\hline
N  & WENO-JS  & WENO-NS  & WENO-P  & MWENO-P\tabularnewline
\hline
10  & 3.1593e-02(---)  & 3.1974e-02(---)  & 3.7323e-02(---)  & 5.2290e-03(---)\tabularnewline
20  & 1.5177e-03(4.37)  & 8.6517e-04(5.27)  & 1.3817e-03(4.75)  & 2.1156e-04(4.62)\tabularnewline
40  & 4.5188e-05(5.06)  & 9.8706e-06(6.45)  & 1.8070e-05(6.25)  & 6.6182e-06(4.99)\tabularnewline
80  & 1.4000e-06(5.01)  & 2.1181e-07(5.54)  & 2.0924e-07(6.43)  & 2.0420e-07(5.01)\tabularnewline
160  & 4.3621e-08(5.00)  & 6.3443e-09(5.06)  & 5.5174e-09(5.24)  & 6.3409e-09(5.00)\tabularnewline
320  & 1.3600e-09(5.00)  & 1.9759e-10(5.00)  & 1.8842e-10(4.87)  & 1.9757e-10(4.99)\tabularnewline
\hline
\end{tabular}\caption{\label{table:2}$L_{1}$ errors of linear advection Eq. (\ref{eq:lintransp})
with initial condition (\ref{48})}
\end{table}

\begin{table}[H]
$\qquad$%
\begin{tabular}{ccccc}
\hline
N  & WENO-JS  & WENO-NS  & WENO-P  & MWENO-P\tabularnewline
\hline
10  & 1.3639e-01(---)  & 7.9335e-02(---)  & 7.8126e-02(---)  & 8.3177e-02(---)\tabularnewline
20  & 1.2790e-02(3.41)  & 1.3137e-02(2.59)  & 1.3174e-02(2.56)  & 5.5774e-03(3.89)\tabularnewline
40  & 1.0952e-03(3.54)  & 3.3350e-04(5.29)  & 3.6552e-04(5.17)  & 2.0964e-04(4.73)\tabularnewline
80  & 8.7557e-05(3.64)  & 6.5660e-06(5.66)  & 7.6657e-06(5.57)  & 6.7007e-06(4.96)\tabularnewline
160  & 7.4148e-06(3.56)  & 2.0946e-07(4.97)  & 1.8501e-07(5.37)  & 2.0988e-07(4.99)\tabularnewline
320  & 4.0271e-07(4.20)  & 6.5526e-09(4.99)  & 6.2591e-09(4.88)  & 6.5526e-09(5.00)\tabularnewline
\hline
\end{tabular}\caption{\label{table:3}$L_{\infty}$ errors of linear advection Eq. (\ref{eq:lintransp})
with initial condition (\ref{49})}
\end{table}

\begin{table}[H]
$\qquad$%
\begin{tabular}{ccccc}
\hline
N  & WENO-JS  & WENO-NS  & WENO-P  & MWENO-P\tabularnewline
\hline
10  & 6.6073e-02(---)  & 2.6674e-02(---)  & 2.6818e-02(---)  & 3.8523e-02(---)\tabularnewline
20  & 4.9673e-03(3.73)  & 3.0373e-03(3.13)  & 3.9893e-03(2.74)  & 2.1165e-03(4.18)\tabularnewline
40  & 3.7068e-04(3.74)  & 7.3314e-05(5.37)  & 7.4578e-05(5.74)  & 7.6161e-05(4.79)\tabularnewline
80  & 1.7135e-05(4.43)  & 2.3943e-06(4.93)  & 1.9211e-06(5.27)  & 2.3652e-06(5.00)\tabularnewline
160  & 7.3448e-07(4.54)  & 7.3574e-08(5.02)  & 6.6582e-08(4.85)  & 7.3534e-08(5.00)\tabularnewline
320  & 2.5137e-08(4.86)  & 2.2922e-09(5.00)  & 2.2006e-09(4.91)  & 2.2927e-09(5.00)\tabularnewline
\hline
\end{tabular}\caption{\label{table:4}$L_{1}$ errors of linear advection Eq. (\ref{eq:lintransp})
with initial condition (\ref{49})}
\end{table}
\subsubsection{Example 2:}
For the partial differential equation (\ref{eq:lintransp}), consider
the initial condition
\begin{equation}
u(x,0)=\text{sin}(\pi x)^{3}.\label{50}
\end{equation}
For this initial condition, the first derivative and second derivatives
vanish whereas the third derivative is non-zero. The $L_{1}$ and
$L_{\infty}$ errors along with the order of convergence for the MWENO-P
scheme is shown in Table \ref{table:5} and Table \ref{table:6} respectively.
The results shows the MWENO-P scheme achieves the desired order of
accuracy where as WENO-NS and WENO-P fails to achieve.
\begin{table}[H]
$\qquad$%
\begin{tabular}{ccccc}
\hline
N  & WENO-JS  & WENO-NS  & WENO-P  & MWENO-P\tabularnewline
\hline
10  & 1.8058e-01(---)  & 1.5556e-01(---)  & 1.4980e-01(---)  & 1.6087e-01(---)\tabularnewline
20  & 6.3274e-02(1.51)  & 3.1225e-02(2.31)  & 3.3257e-02(2.17)  & 2.6674e-02(2.59)\tabularnewline
40  & 6.0354e-03(3.39)  & 6.6848e-03(2.22)  & 8.1907e-03(2.02)  & 3.0858e-03(3.11)\tabularnewline
80  & 9.1031e-04(2.72)  & 9.1944e-04(2.86)  & 1.2103e-03(2.75)  & 3.0617e-04(3.33)\tabularnewline
160  & 4.8182e-05(4.23)  & 1.0142e-04(3.18)  & 1.3156e-04(3.20)  & 1.3278e-06(7.84)\tabularnewline
320  & 8.0849e-07(4.89)  & 9.8843e-06(3.35)  & 1.3798e-05(3.25)  & 3.8780e-08(5.09)\tabularnewline
640  & 1.3257e-08(5.93)  & 9.9097e-07(3.31)  & 1.3735e-06(3.32)  & 1.1322e-09(5.09)\tabularnewline
1280  & 2.3166e-10(5.83)  & 9.9188e-08(3.32)  & 1.4340e-07(3.26)  & 3.4832e-11(5.02)\tabularnewline
\hline
\end{tabular}\caption{\label{table:5}$L_{1}$ errors of linear advection Eq. (\ref{eq:lintransp})
with initial condition (\ref{50})}
\end{table}

\begin{table}[H]
$\qquad$%
\begin{tabular}{ccccc}
\hline
N  & WENO-JS  & WENO-NS  & WENO-P  & MWENO-P\tabularnewline
\hline
10  & 2.4666e-01(---)  & 2.2154e-01(---)  & 2.1397e-01(---)  & 2.2637e-01(---)\tabularnewline
20  & 1.2969e-01(0.92)  & 6.0015e-02(1.88)  & 7.2114e-02(1.56)  & 4.7482e-02(2.25)\tabularnewline
40  & 1.2970e-02(3.32)  & 1.6773e-02(1.83)  & 2.1148e-02(1.77)  & 1.1456e-02(2.05)\tabularnewline
80  & 3.8067e-03(1.76)  & 3.7616e-03(2.15)  & 5.3389e-03(1.98)  & 2.2518e-03(2.34)\tabularnewline
160  & 3.4089e-04(3.48)  & 8.0170e-04(2.23)  & 1.0032e-03(2.41)  & 4.2414e-06(9.05)\tabularnewline
320  & 6.9825e-06(5.60)  & 1.4264e-04(2.49)  & 1.9960e-04(2.32)  & 9.6777e-08(5.45)\tabularnewline
640  & 7.5668e-08(6.52)  & 2.4583e-05(2.53)  & 3.6812e-05(2.43)  & 1.7514e-09(5.78)\tabularnewline
1280  & 7.2140e-10(6.71)  & 4.7171e-06(2.38)  & 6.7828e-06(2.44)  & 5.4793e-11(4.99)\tabularnewline
\hline
\end{tabular}\caption{\label{table:6}$L_{\infty}$ errors of linear advection Eq. (\ref{eq:lintransp})
with initial condition (\ref{50})}
\end{table}

\subsubsection{Example 3:}
\begin{figure}[H]
\centering{}\includegraphics[width=18cm,height=7cm]{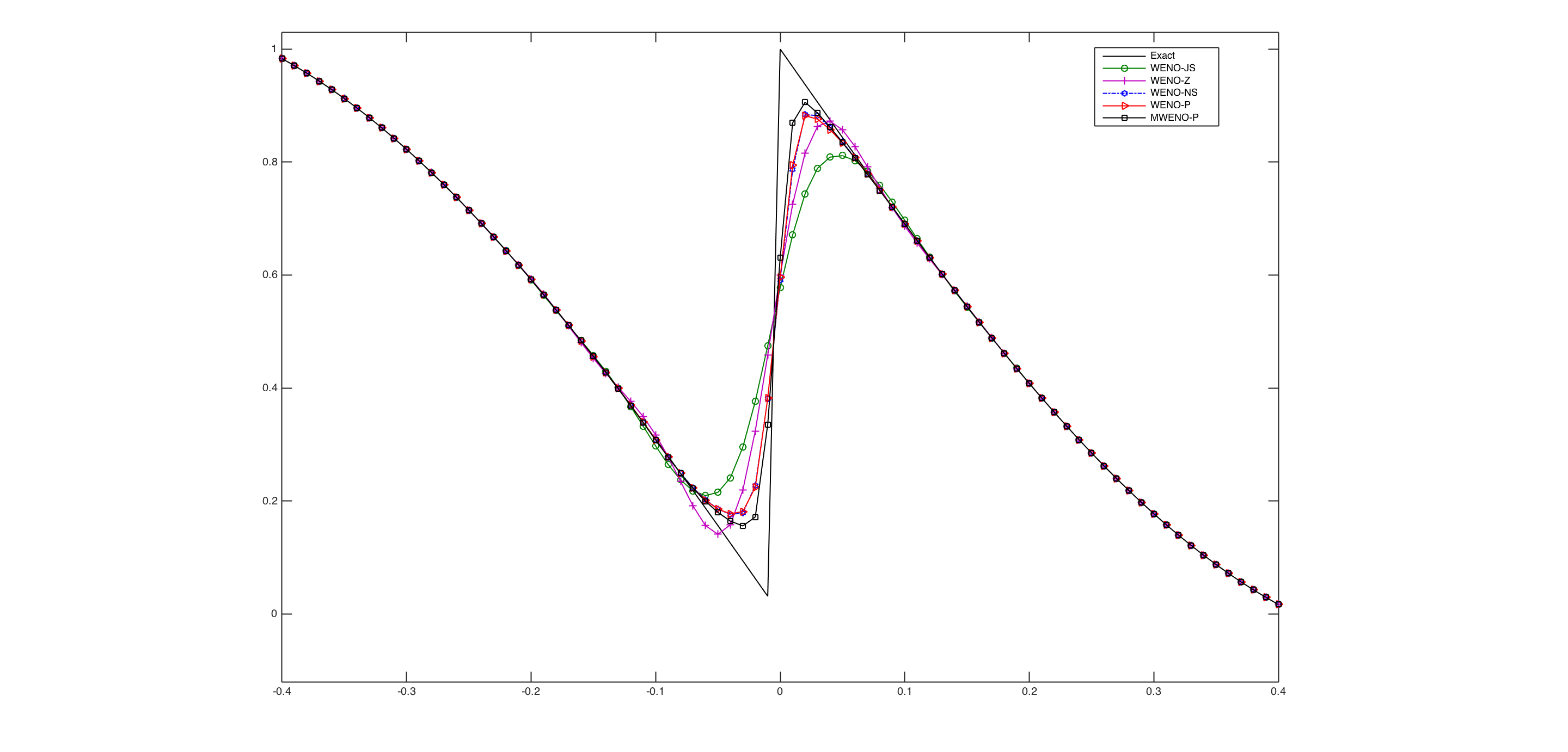}\caption{\label{fig:sin3x}Numerical solution of linear advection Eq. (\ref{eq:lintransp})
with the initial condition (\ref{52}) }
\end{figure}

For linear advection equation (\ref{eq:lintransp}), let the initial
condition be
\begin{equation}
u(x,0)=u_{0}(x)=\begin{cases}
-\text{sin}(\pi x)-\frac{1}{2}x^{3} & \text{ for }\,\,-1\leq x<0,\\
-\text{sin}(\pi x)-\frac{1}{2}x^{3}+1 & \text{ for }\,\,0\leq x<1.
\end{cases}\label{52}
\end{equation}
which is a piecewise sine function with jump discontinuity at $x=0.$
The solution is computed with the CFL number $0.5$ with uniform concretization
of the domain and the step size is $\triangle x=0.01$ up to time
$t=8.$ The approximate solution computed with MWENO-P along with
WENO-JS, WENO-NS and WENO-P schemes is plotted in Figure \ref{fig:sin3x}
against the exact solution. It can be observed from the plot that
the proposed scheme performs better than other schemes near the jump
discontinuity.
\subsubsection{Example 4:}
\begin{figure}[H]
\centering{}\includegraphics[width=18cm,height=7cm]{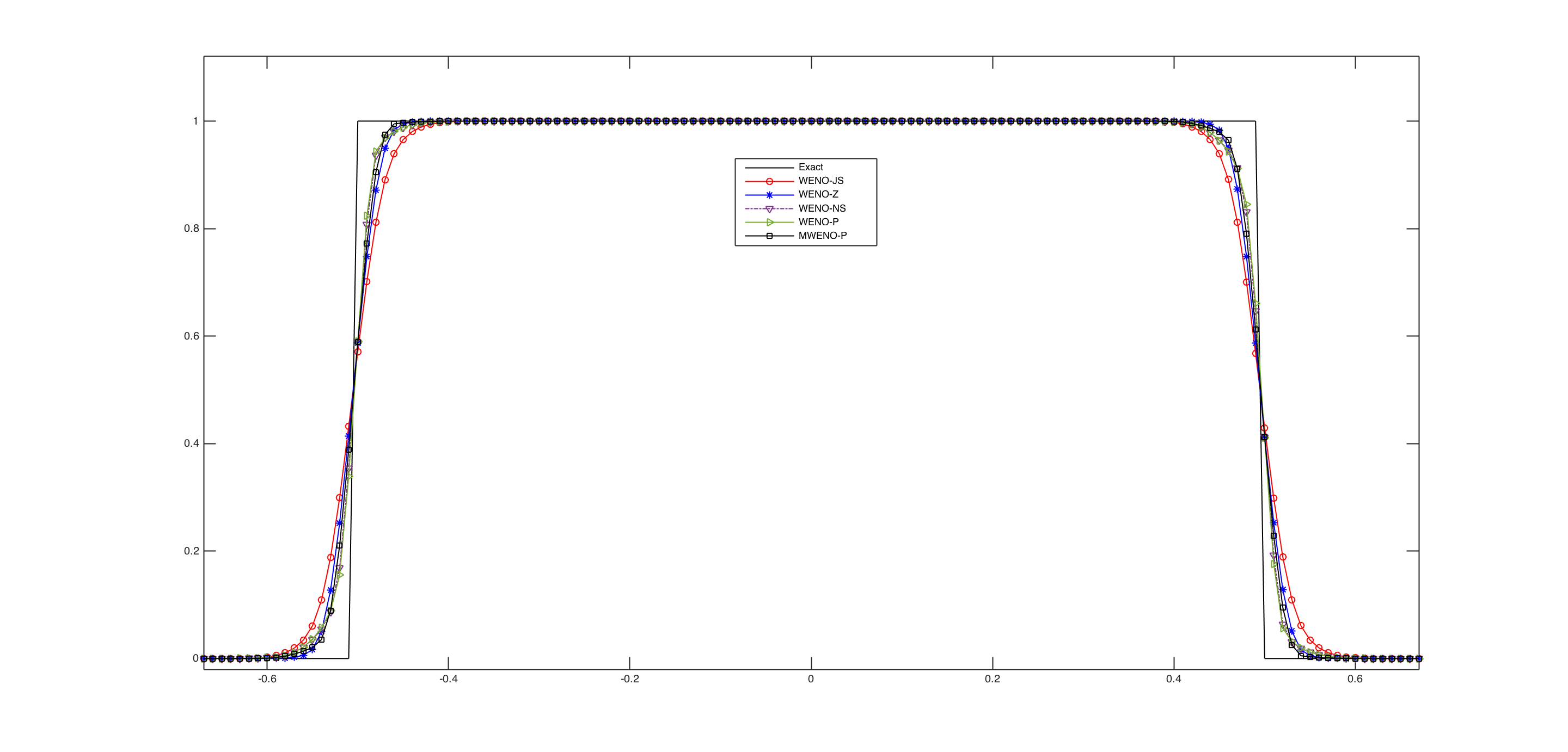}\caption{\label{discontIC}Numerical solution of linear advection Eq. (\ref{eq:lintransp})
with the initial condition (\ref{54}) }
\end{figure}

Consider the discontinuous profile
\begin{equation}
u(x,0)=\begin{cases}
1 & \text{ if }-0.5\leq x<0.5,\\
0 & \text{ otherwise. }
\end{cases}\label{54}
\end{equation}
for the linear equation (\ref{eq:lintransp}). The equation is solved
with the CFL number $0.5$ with the spatial step size $\triangle x=0.01$.
For time $t=10$ the computed approximate solutions are plotted against
exact solution in Figure \ref{discontIC}. The proposed scheme MWENO-P
has better approximation than WENO-JS, WENO-NS and WENO-P schemes
especially near the areas of discontinuities.
\subsubsection{Example 5:}
\begin{figure}[H]
\centering{}\includegraphics[width=18cm,height=8cm]{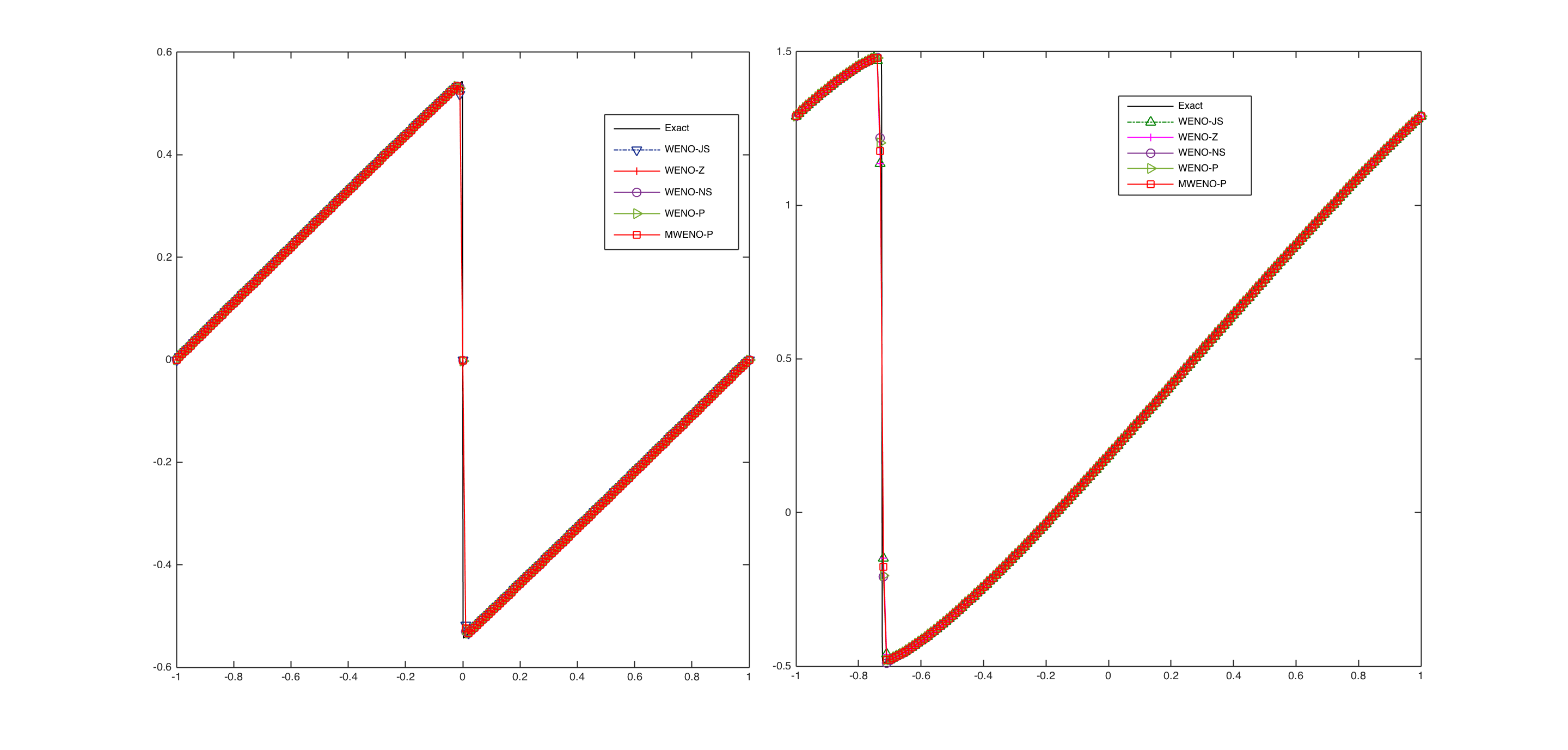}\caption{\label{fig:3-1}Approximate solution of (\ref{eq:48}) with initial
condition (\ref{56}) and (\ref{57}) }
\end{figure}

The WENO schemes are designed as a shock capturing schemes for solving
the hyperbolic conservation laws, for this as a last example in this
section considered the Burger's equation
\begin{align}
 & u_{t}+(\frac{u^{2}}{2})_{x}=0,\,\,-1\leq x\leq1,\,\,t\geq0,\label{eq:48}\\
 & u(x,0)=u_{0}(x).\nonumber
\end{align}
subject to periodic boundary conditions. In Figure \ref{fig:3-1}, the numerical result of the fifth-order WENO schemes for the initial
conditions
\begin{eqnarray}
 & u_{0}(x)= & -\text{sin}(\pi x).\label{56}
\end{eqnarray}
at time $t=1.5$ and
\begin{eqnarray}
 & u_{0}(x)= & \frac{1}{2}+\text{sin}(\pi x).\label{57}
\end{eqnarray}
at $t=0.55$ are plotted respectively. The exact or reference solution
is calculated with $2000$ grid points with WENO-JS scheme and the
approximate solutions are computed with $200$ grid points in space.
It is shown that the shocks are very well captured by all the schemes.

\subsection{Euler equations in one space dimension}

\hspace{1cm}The one-dimensional Euler equations are given by
\begin{equation}
\begin{pmatrix}\rho\\
\rho u\\
E
\end{pmatrix}_{t}+\begin{pmatrix}\rho u\\
\rho u^{2}+p\\
u(E+p)
\end{pmatrix}_{x}=0\label{eq:60}
\end{equation}
where $\rho,u,E,p$ are the density, velocity, total energy and pressure
respectively. The system (\ref{eq:60}) represents the conservation
of mass, momentum and energy. The total energy for an ideal polytropic
gas is defined as
\[
E=\frac{p}{\gamma-1}+\frac{1}{2}\rho u^{2},
\]
where $\gamma$ is the ratio of specific heats and its value is taken
as $\gamma=1.4.$

The numerical flux $\hat{f}_{j+\frac{1}{2}}$ at $x_{j+\frac{1}{2}}$
is calculated based on the steps given in \cite{shunotes}, which
are reproduced here for completeness.
\begin{itemize}
\item Compute an average state $u_{j+\frac{1}{2}}$ using Roe's mean matrix.
\item Compute the left eigenvectors, the right eigenvectors and the eigenvalues
of the Jacobian $f^{'}(u)$ at the average state $u_{j+\frac{1}{2}}$
and denote them by $L=L(u_{j+\frac{1}{2}})$, $R=R(u_{j+\frac{1}{2}}),$
and $\wedge=\wedge(u_{j+\frac{1}{2}})$ respectively. Note that $L=R^{-1}$.
\item Projecting the flux and the solution which are in the potential stencil
of the WENO reconstruction for obtaining the flux $\hat{f}_{j+\frac{1}{2}}$
to the local characteristic fields by $s_{j}=R^{-1}u_{j}$, $q_{j}=R^{-1}f(u_{j})$,
$j$ in a neighborhood of $i$.
\item To obtain the corresponding component of the flux $\hat{q}_{j+\frac{1}{2}}^{\pm},$
Lax-Friedrich's flux splitting and the WENO reconstruction procedure
is used for each component of the characteristic variables.
\item Projecting back the characteristic variables to physical variables
by $\hat{f}_{j+\frac{1}{2}}^{\pm}=R\hat{q}_{j+\frac{1}{2}}^{\pm}$.
\item Finally, form the flux by taking $\hat{f}_{j+\frac{1}{2}}=\hat{f}_{j+\frac{1}{2}}^{+}+\hat{f}_{j+\frac{1}{2}}^{-}$.
\end{itemize}
For time evaluation, we used third-order TVD Runge-Kutta scheme (\ref{eq:45}).
Now consider the one dimensional Riemann problem for Euler system
of equations (\ref{eq:60}) i.e., with the initial condition
\[
U(x,0)=\begin{cases}
U_{L} & \text{ if }x<x_{0},\\
U_{R} & \text{ if }x\geq x_{0}.
\end{cases}
\]
where $U_{L}=(\rho_{l},u_{l},p_{l})$ and $U_{R}=(\rho_{r},u_{r},p_{r})$.
In the following, various test cases are taken for numerical study
of MWENO-P scheme.

\subsubsection{Sod's shock tube problem: }

\begin{figure}[!ht]
\centering{}\includegraphics[width=18cm,height=8cm]{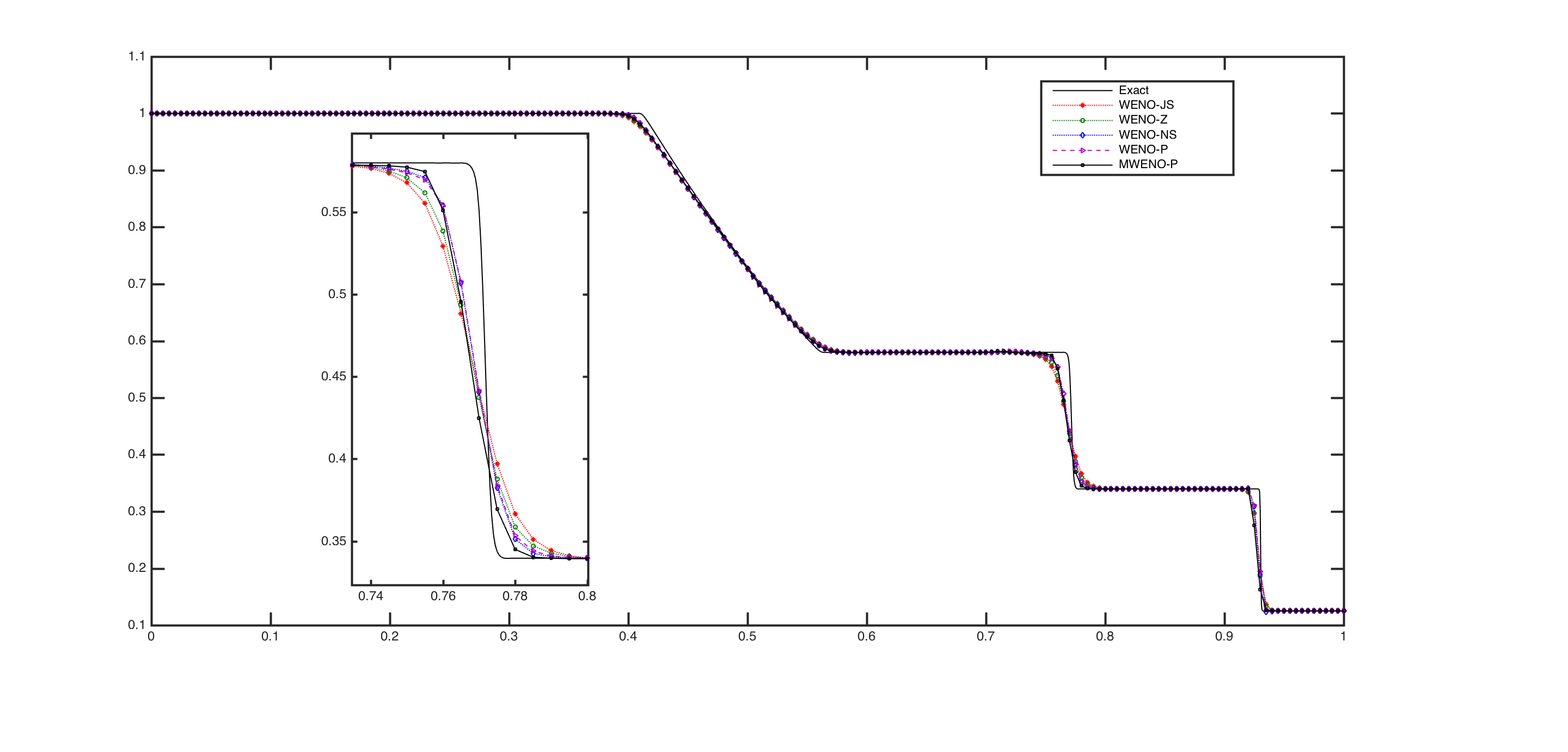}\caption{\label{fig:Den-sod1}Density profile for Sod's problem and the zoomed
region at the contact discontinuity}
\end{figure}

\begin{figure}[!ht]
\centering{}\includegraphics[width=18cm,height=8cm]{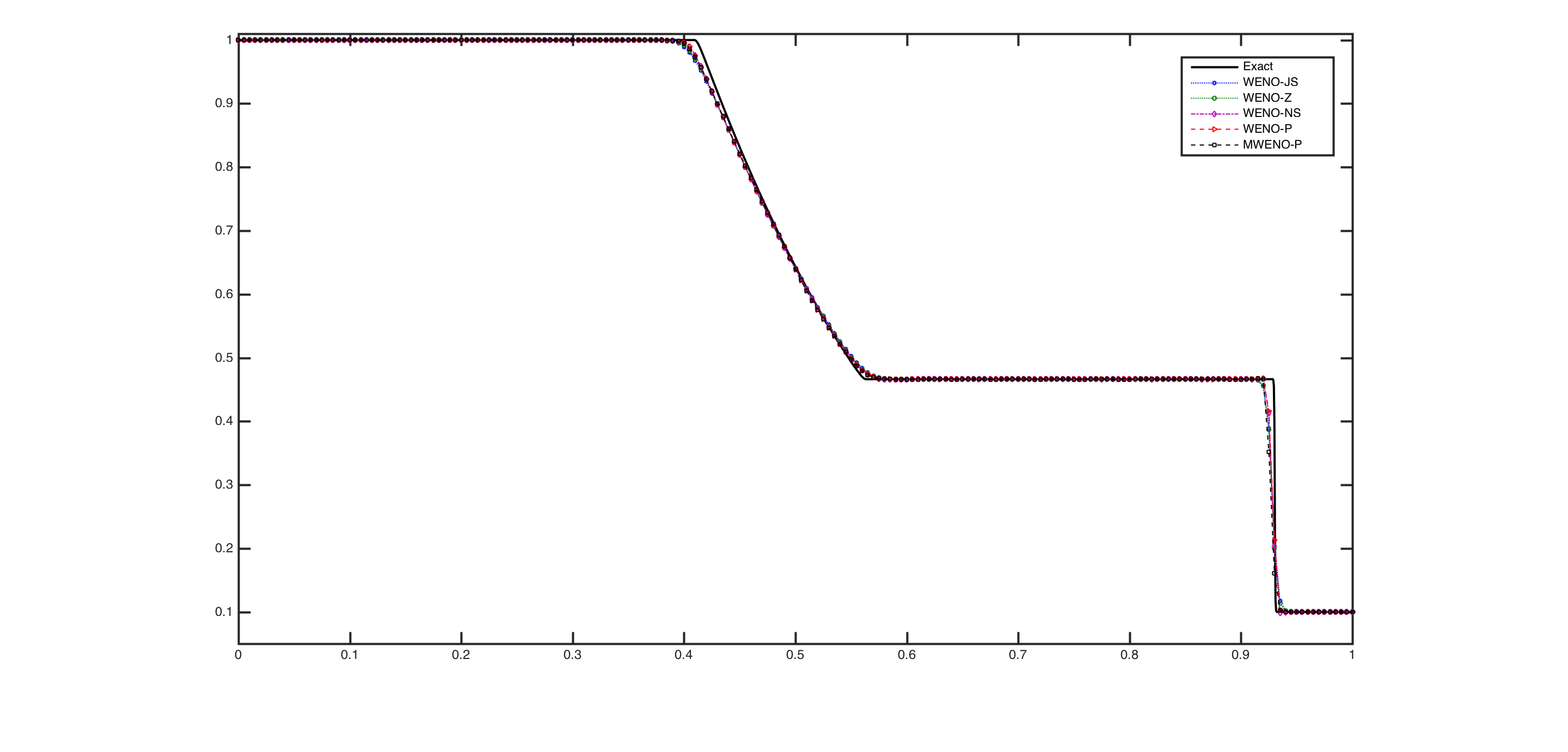}\caption{\label{fig:Pres-sod1}Pressure profile of Sod's shock tube problem }
\end{figure}

For this the initial condition is given by
\[
U(x,0)=\begin{cases}
(1.0,0.75,1.0), & \text{ if }0\leq x\leq0.5,\\
(0.125,0.0,0.1), & \text{ if }0.5\leq x\leq1.
\end{cases}
\]
This is an example of modified version of Sod's problem defined in
\cite{Sod} and its solution contains a right shock wave, a right
traveling contact wave and a left sonic rarefaction wave. Transmissive
boundary conditions are taken for numerical evaluation. The solution
is computed up to time $t=0.2$ with $200$ grid points in space with
the CFL number $0.5.$ The density and pressure profile for various
fifth-order WENO schemes are shown in Figure \ref{fig:Den-sod1} and
Figure \ref{fig:Pres-sod1} respectively against the reference solution
which is calculated with $2000$ grid points using WENO-JS scheme.
It is observed that the proposed scheme MWENO-P performs better than
other WENO schemes near the region of contact discontinuity.

\subsubsection{Shock tube problem with Lax initial condition:}

\begin{figure}[H]
\centering{}\includegraphics[width=18cm,height=8cm]{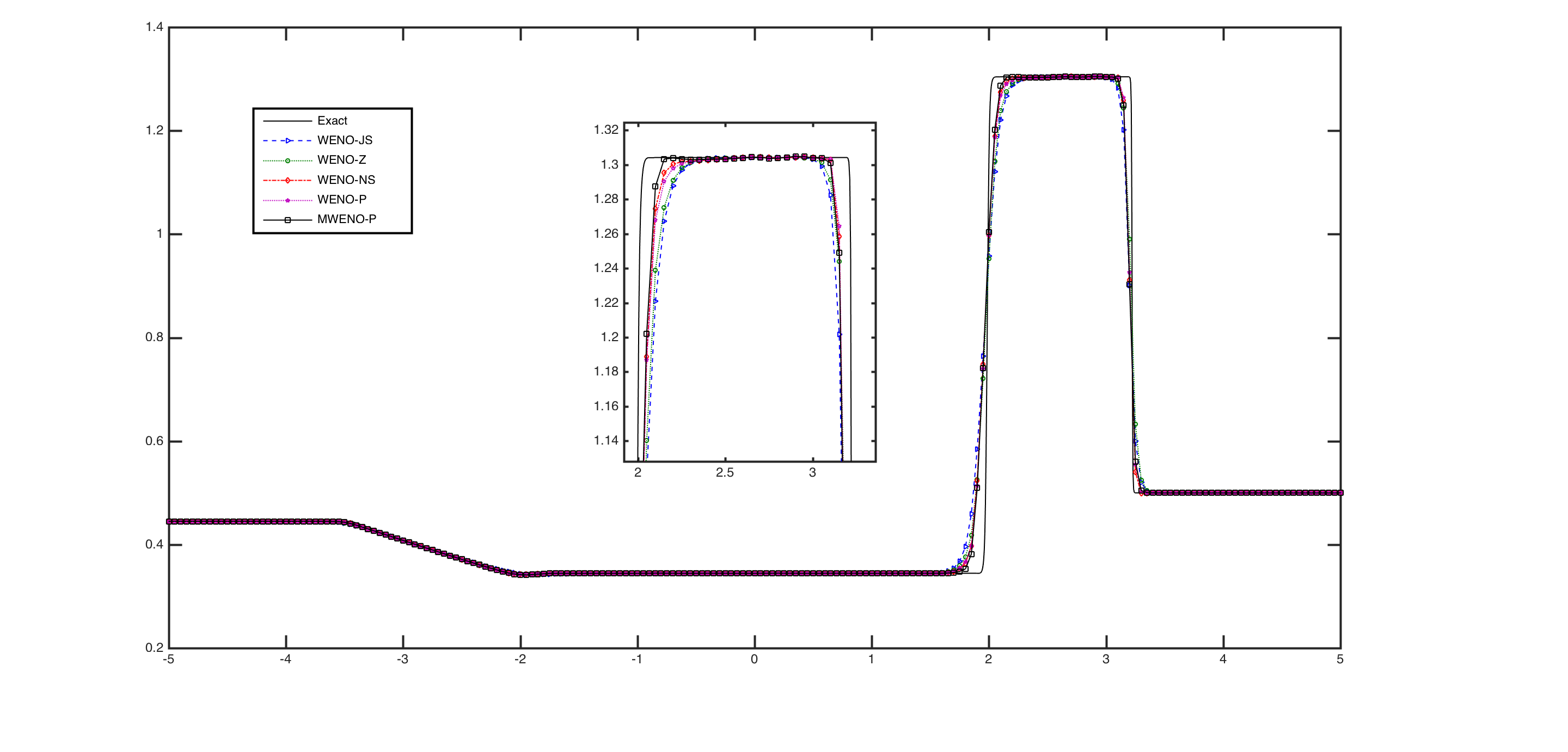}\caption{\label{fig:Den-Lax}Density profile for Lax initial condition}
\end{figure}

\begin{figure}[H]
\centering{}\includegraphics[width=18cm,height=8cm]{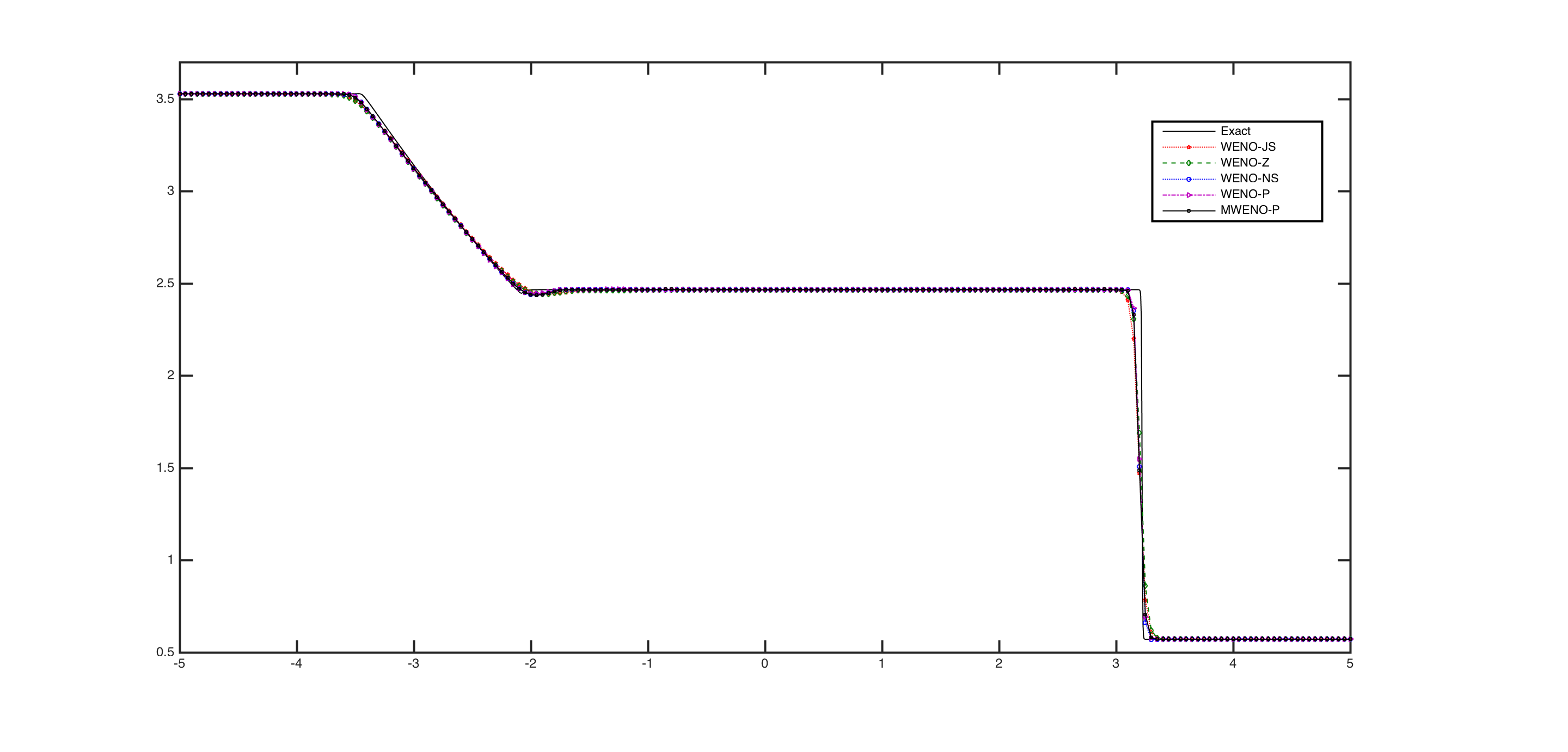}\caption{\label{fig:Pres-Lax}Pressure profile for Lax initial condition}
\end{figure}

For this the initial condition which is given by
\[
U(x,0)=\begin{cases}
(0.445,0.698,3.528), & \text{ if }-5\leq x<0,\\
(0.500,0.000,0.571), & \text{ if }0\leq x\leq5.
\end{cases}
\]
For this, the solution is computed up to time $t=1.3$ with $200$
grid points along space direction and the CFL number is set to be
$0.5$. The reference solution is calculated with $2000$ grid points
by using WENO-JS scheme. For assessing the performance of MWENO-P
scheme, the density profile is plotted against the reference solution
along with other schemes in Figure \ref{fig:Den-Lax}, it can be seen
that MWENO-P scheme performs better than other schemes near the discontinuities.
The pressure profile is displayed in Figure \ref{fig:Pres-Lax}.

\subsubsection{1D shock entropy wave interaction problem: }

The initial condition as given in \cite{woodward} for shock entropy
wave interaction problem is
\[
U(x,0)=\begin{cases}
(3.857143,2.629369,10.33333), & \text{ if }-5\leq x<-4,\\
(1+\epsilon sin(kx),0.000,1.000), & \text{ if }-4\leq x\leq5.
\end{cases}
\]
where $\epsilon$ and $k$ are the amplitude and wavenumber of the
entropy wave respectively, chosen $\epsilon=0.2$ and $k=5.$ This
problem has a right-moving supersonic (Mach 3) shock wave which interacts
with sine waves in a density disturbance, that generates a flow field
with both smooth structures and discontinuities. This flow induces
wave trails behind a right-going shock with wave numbers higher than
the initial density-variation wavenumber $k$. Since the exact solution
is unknown, the reference solution with zero gradient boundary conditions
is obtained by using the fifth-order WENO-JS scheme with $2000$ grid
points. The initial condition contains a jump discontinuity at $x=-4$
and especially, the initial density profile has oscillations on $[-4,5]$.
The solution is computed up to time $t=1.8$ with $200$ spatial grid
points and the CFL number is set to be $0.5$. The Figure $8$ contains
the graph of the solutions calculated at $200$ grid points and its
zoomed region near oscillations. It is observed that MWENO-P scheme
performs better than WENO-JS and WENO-P at $200$ grid points.

\begin{figure}[H]
\centering{}\includegraphics[width=18cm,height=8cm]{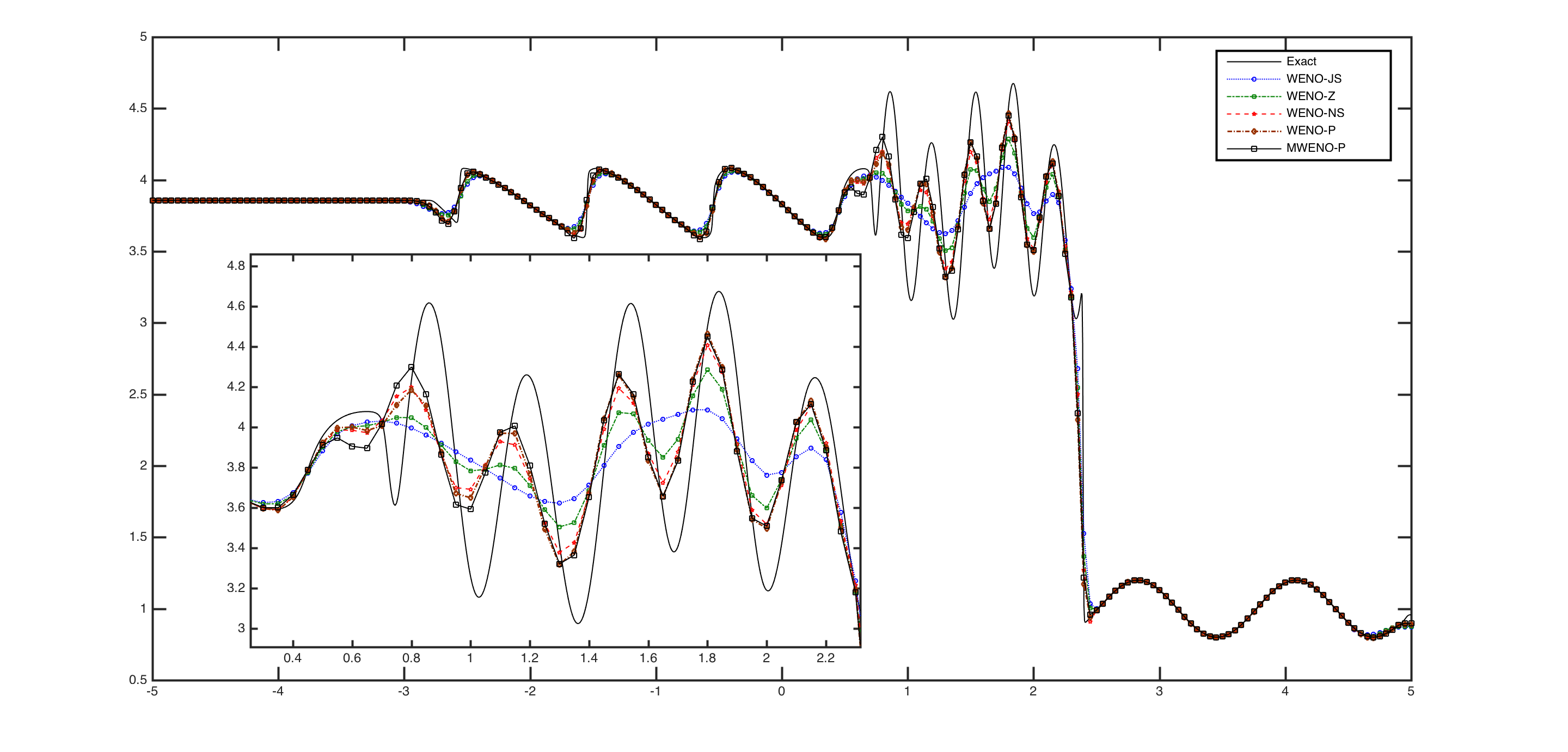}\caption{Shock entropy wave interaction test with 200 grid points}
\end{figure}

\subsection{Two-Dimensional Euler equations}

\hspace{1cm}In this section, we apply the proposed scheme to two-dimensional
problem in cartesian coordinates. The governing two-dimensional compressible
Euler equations is given by

\[
U_{t}+F(U)_{x}+G(U)_{y}=0
\]
where $U=(\rho,\rho u,\rho v,E)^{T}$, $F(U)=(\rho u,P+\rho u^{2},\rho uv,u(E+P))^{T}$,
$G(U)=(\rho v,\rho uv,P+\rho v^{2},v(E+P))^{T}$. The total energy
$E$ and the pressure $p$ is defined by

\[
p=(\gamma-1)(E-\frac{1}{2}\rho(u^{2}+v^{2}))
\]
where $\gamma$ is the ratio of specific heats. Here $\rho,u,v$ are
density, $x$-wise-velocity component and $y$-wise-velocity component
respectively.

\subsubsection{2D problem of 2D gas dynamics}

\begin{figure}[H]
\centering{}\includegraphics[width=18cm,height=8cm]{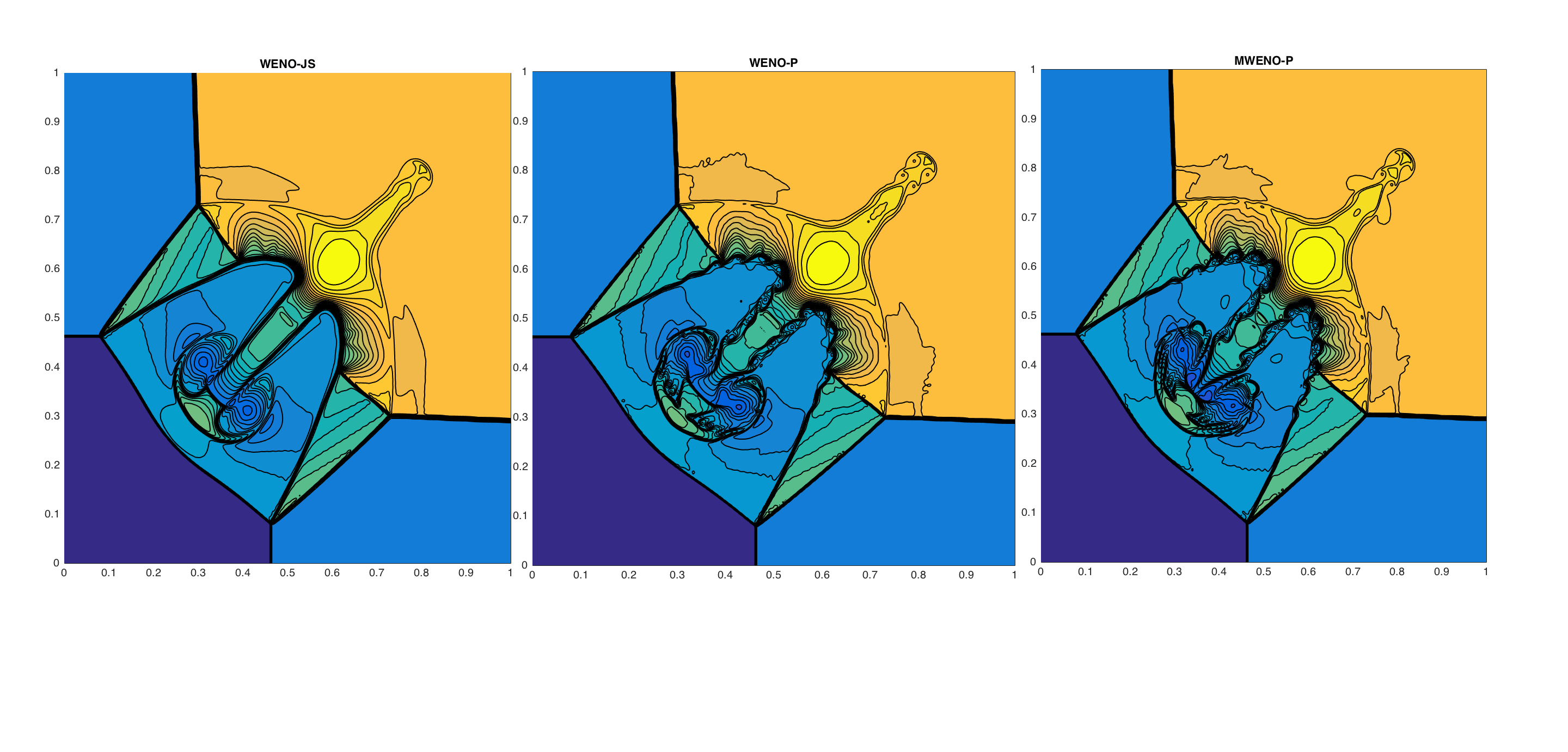}\caption{Density profile of 2D Riemann problem of 2D gas dynamics with the mesh $\Delta x=\Delta y=1/400$}\label{fig:2Dgas}
\end{figure}

\hspace{1cm}The 2D Riemann problem of gas dynamics is proposed in
\cite{Rinne} and solved on the rectangular domain $[0,1]\times[0,1]$.
The 2D Riemann problem is defined by initial constant states in each
quadrant which is divided by lines $x=0.8$ and $y=0.8$ on the square
as
\begin{eqnarray*}
(\rho,u,v,p) & =\begin{cases}
\begin{array}{lr}
(1.5,0,0,1.5) & if\;0.8\leq x\leq1,0.8\leq y\leq1,\\
(0.5323,1.206,0,0.3) & if\;0\leq x<0.8,0.8\leq y\leq1,\\
(0.138,1.206,1.206,0.029) & if\;0\leq x<0.8,0\leq y<0.8,\\
(0.5323,0,1.206,0.3) & if\;0.8<x\leq1,0\leq y\leq0.8,
\end{array}\end{cases}
\end{eqnarray*}
with Dirichlet boundary conditions. According to the initial conditions,
four shocks come into being and produce a narrow jet. The numerical
solution is calculated with $400\times400$ grid points up to $t=0.8$
with the CFL number $0.5$. An examination of these results reveals in figure \ref{fig:2Dgas}
that, MWENO-P yields a better solution of the complex structure appearing
when compares to WENO-JS and WENO-P schemes.

\subsubsection{Two-Dimensional Rayleigh-Taylor instability}

\hspace{1cm}The Taylor instability happens on an interface between
two fluids of different densities when an acceleration is directed
from heavier fluid to lighter fluid. This problem has been simulated
extensively in the literature \cite{young}. The computational domain
is $[0,1/4]\times[0,1]$ and the initial conditions are
\[
(\rho,u,v,p)=\begin{cases}
\begin{array}{cc}
(2,0,-0.025acos(8\pi x),2y+1) & \text{if }0\leq y<0.5,\\
(1,0,-0.025acos(8\pi x),y+\frac{3}{2}) & \text{if }0.5\leq y<1.
\end{array}\end{cases}
\]
with the sound speed $a=\sqrt{\gamma p/\rho}$ and the ratio of specific
heats $\gamma=5/3$. The gravitational effect is introduced by adding
$\rho$ and $\rho v$ to the right hand side of third equation and
fourth equation, respectively of the two-dimensional Euler equation.
Reflective boundary conditions are imposed for the left and right
boundaries and

\[
(\rho,u,v,p)=\begin{cases}
\begin{array}{cc}
(2,0,0,2.5) & \text{top boundary,}\\
(1,0,0,1) & \text{bottom boundary.}
\end{array}\end{cases}
\]
The final simulation time is $t=1.95$. The density contour plotted
in Fig. 10 shows that the WENO-P and MWENO-P obtains more complex
structures than the other schemes.

\begin{figure}[H]
\centering{}\includegraphics[width=16cm,height=10cm]{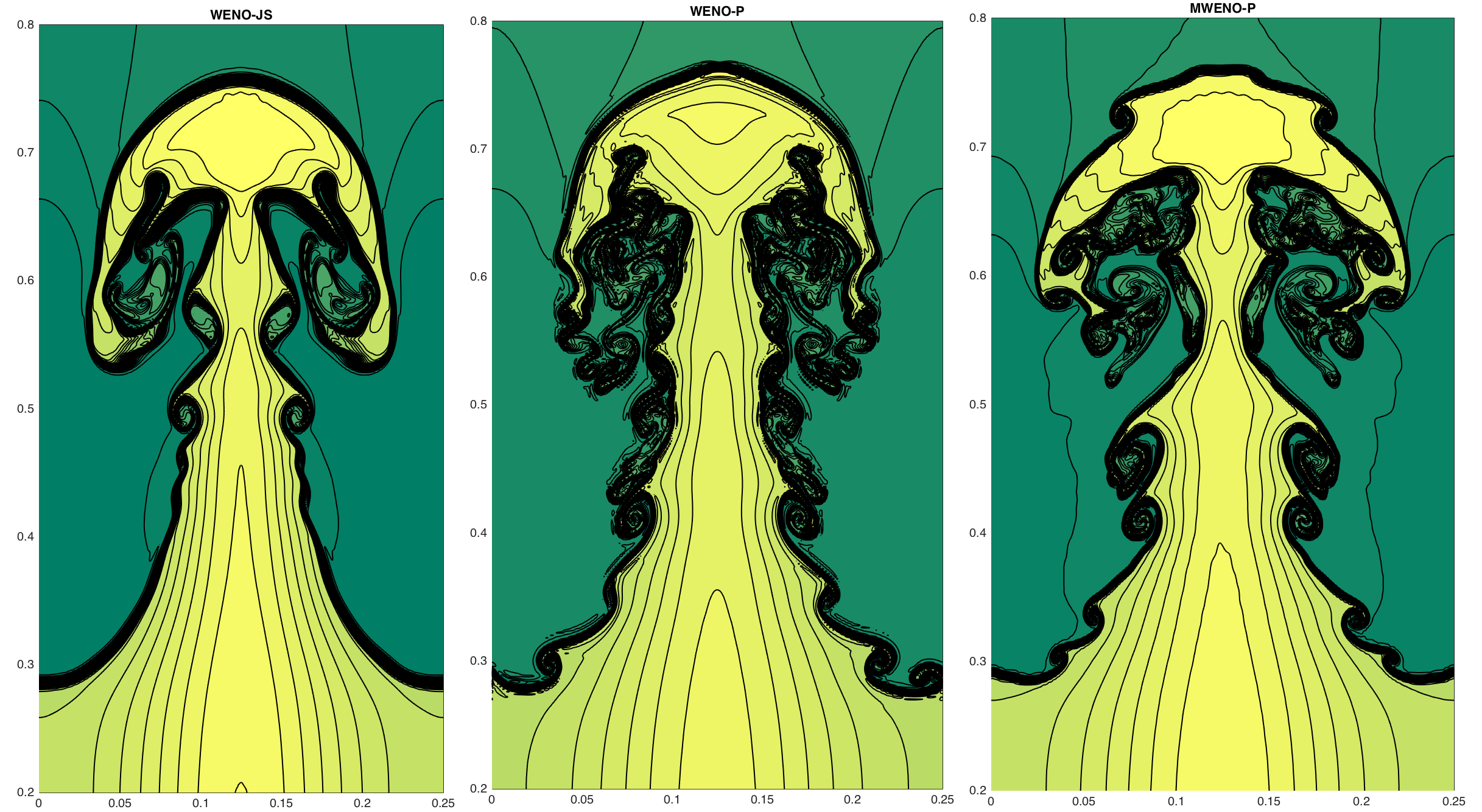}\protect\caption{Density profiles of Rayleigh-Taylor instability problem at t=1.95 with
the mesh $\Delta x=\Delta y=1/500$}
\end{figure}

\subsubsection{Double Mach reflection of a strong shock}

\hspace{0.6cm}The final test considered in this paper is the two-dimensional
double mach reflection problem\cite{woodward} where a vertical shock
wave moves horizontally into a wedge that is inclined by some angle.
The computational domain for this problem is chosen to be $[0,4]\times[0,1],$
and the reflecting wall lies at the bottom of the computational domain
for $\frac{1}{6}\leq x\leq4$. Initially a right-moving Mach $10$
shock is positioned at $x=\frac{1}{6},y=0,$ and makes a $60^{o}$
angle with the x-axis. For the bottom boundary, the exact postshock
condition is imposed for the part from $x=0$ to $x=\frac{1}{6}$
and a reflective boundary condition is used for the rest. The top
boundary of our computational domain uses the exact motion of the
Mach $10$ shock. Inflow and outflow boundary conditions are used
for the left and right boundaries. The unshocked fluid has a density
of $1.4$ and a pressure of $1.$ The problem was run till $t=0.2$
and the blow-up region around the double Mach stems. The ration of
specific heats $\gamma=1.4$ and we set CFL number as $0.5.$ The
results in $[0,3]\times[0,1]$ are displayed for the WENO-JS, WENO-P
and MWENO-P schemes in figures 11, 12 and 13 respectively. The figure
14 shows that the performance of the schemes WENO-JS, WENO-P and MWENO-P
at the Mach stem of the density variable at the final time with several
grid points. It can be clearly seen that MWENO-P resolves better the
instabilities around the Mach stem.

\begin{figure}[H]
\centering{}\includegraphics[width=16cm,height=6cm]{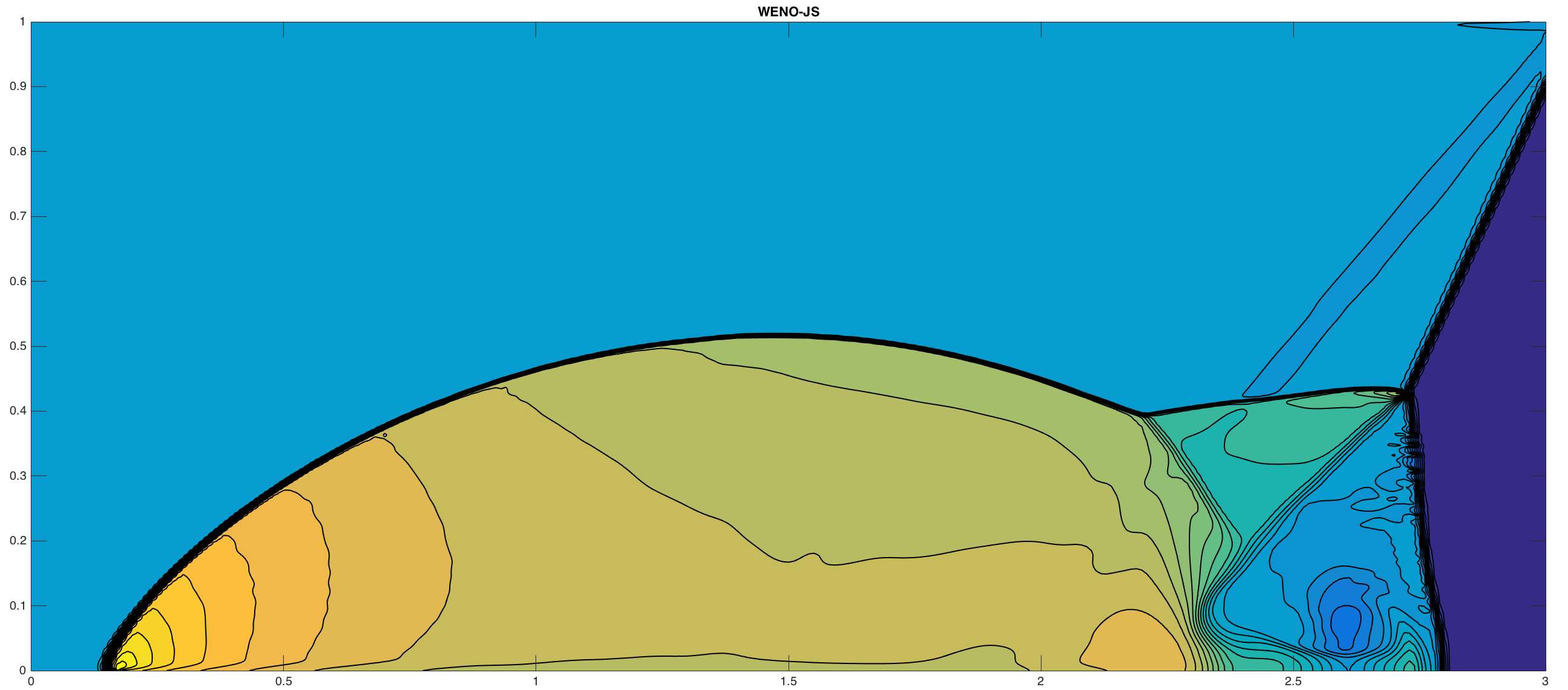}\protect\caption{Density profiles of Double Mach reflection of a strong shock with
WENO-JS scheme at t=0.2 with the mesh $\Delta x=\Delta y=1/400$}
\end{figure}
\begin{figure}[H]
\centering{}\includegraphics[width=16cm,height=6cm]{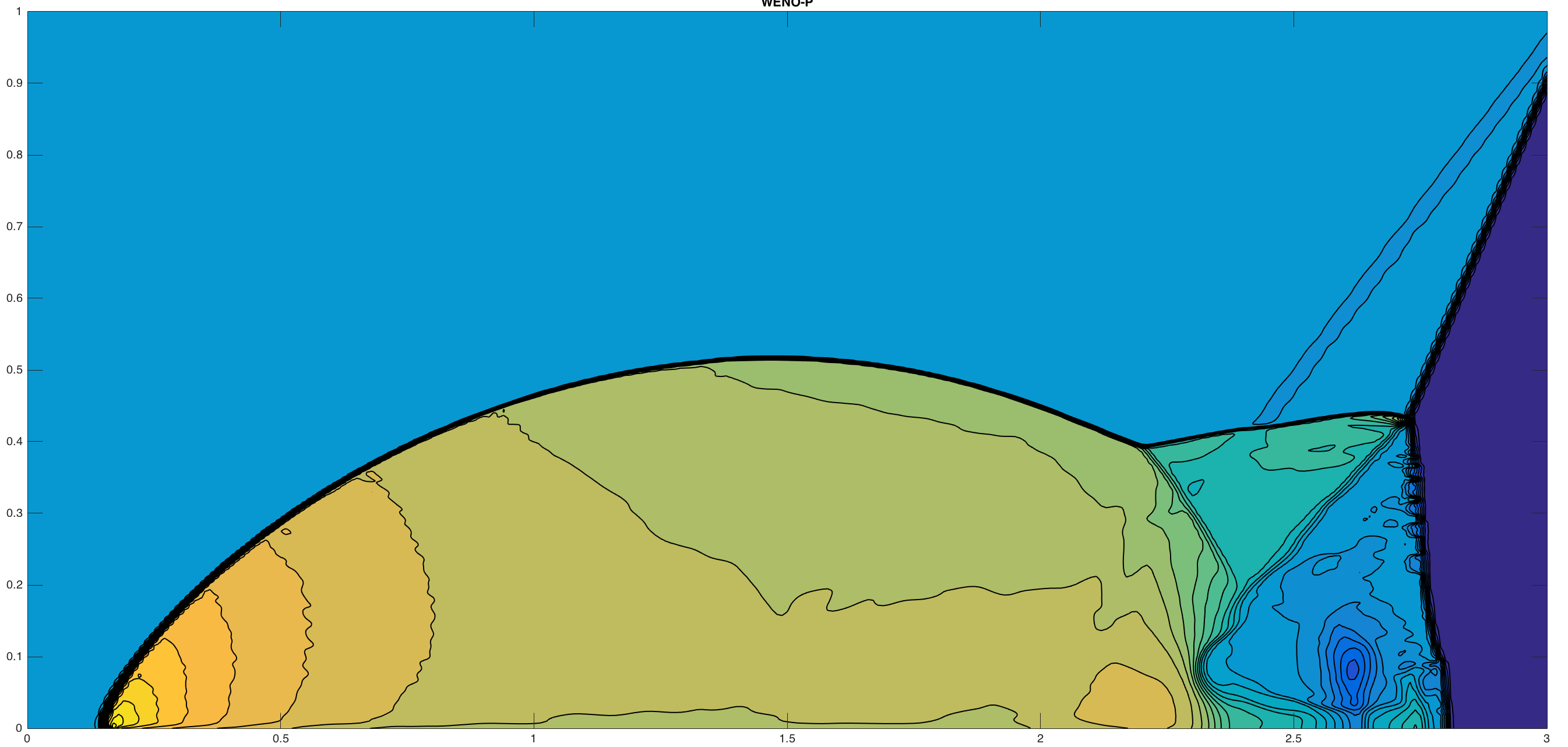}\protect\caption{Density profiles of Double Mach reflection of a strong shock with
WENO-P scheme at t=0.2 with the mesh $\Delta x=\Delta y=1/400$}
\end{figure}
\begin{figure}[H]
\centering{}\includegraphics[width=16cm,height=6cm]{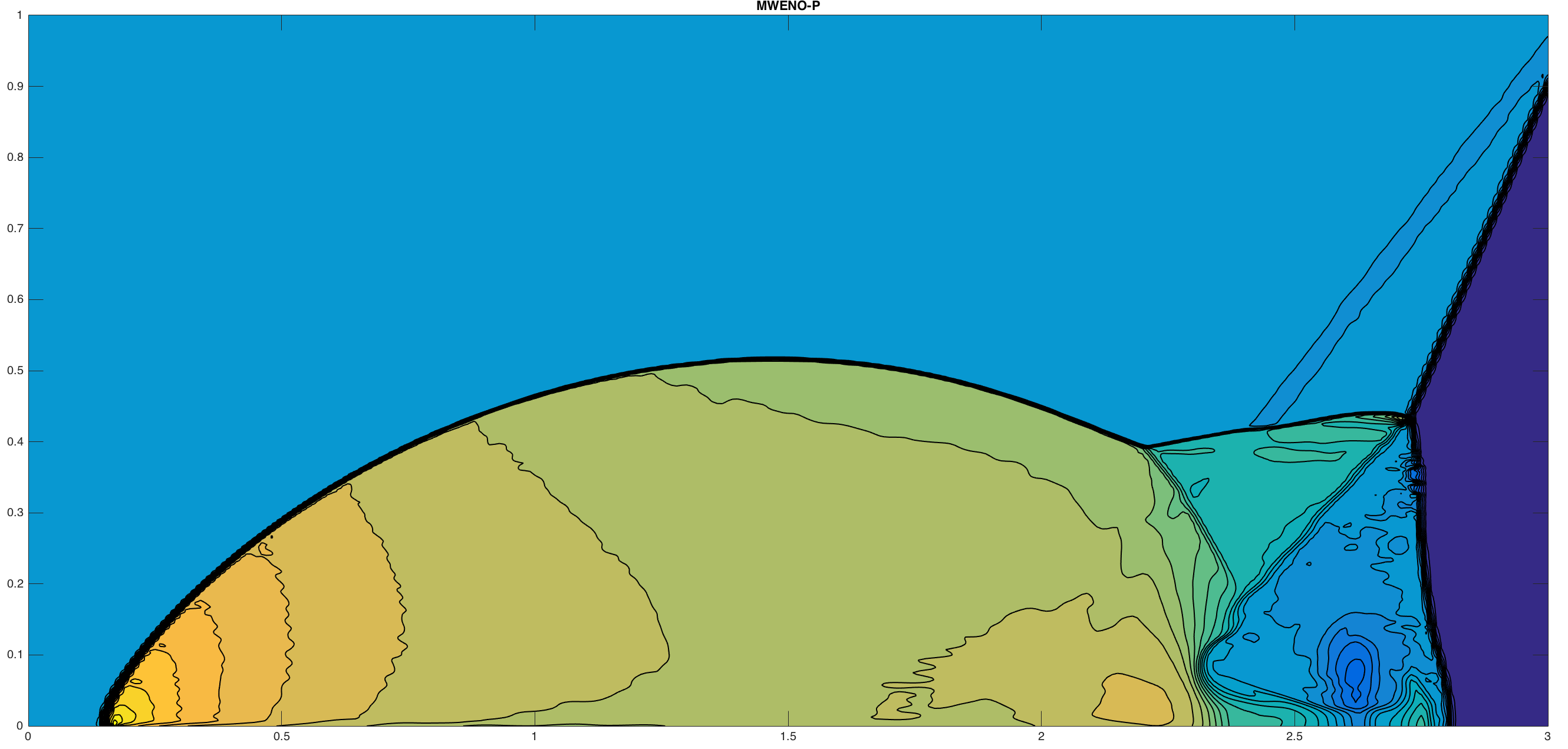}\protect\caption{Density profiles of Double Mach reflection of a strong shock with
MWENO-P scheme at t=0.2 with the mesh $\Delta x=\Delta y=1/400$}
\end{figure}
\begin{figure}[H]
\centering{}\includegraphics[width=16cm,height=8cm]{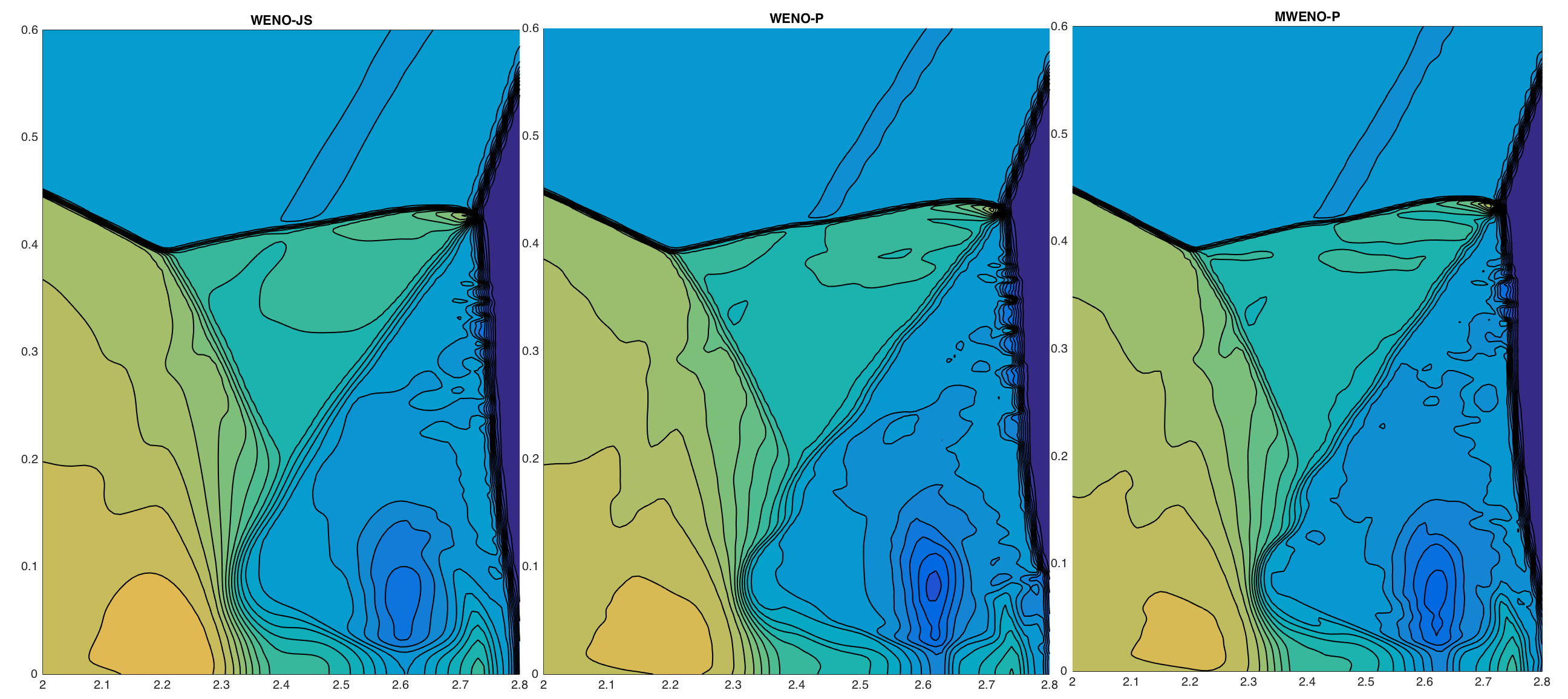}\protect\caption{Density profiles of Double Mach reflection of a strong shock at Mach
stem with WENO-JS, WENO-P and MWENO-P schemes at t=0.2 with the mesh
$\Delta x=\Delta y=1/400$}
\end{figure}

\section{Conclusion}

\hspace{0.6cm}In this paper, a modified fifth-order WENO scheme to
approximate the solution of nonlinear hyperbolic conservation laws
named as MWENO-P was presented, by introducing a new global smoothness
indicator based on undivided differences of second order derivatives
of an interpolation polynomial over a stencil. The motivation to the
study is that the WENO-NS and WENO-P schemes do not satisfy the sufficient
condition if first-order derivatives vanish but the second derivative
is non-zero. The proposed scheme satisfied the sufficient condition
even if first and second derivatives vanish but the third derivative
is non-zero. The approximate solutions to the one-dimensional scalar,
system and two-dimensional system of hyperbolic conservation laws
are simulated with the proposed scheme and compared it with other
fifth-order WENO schemes. Numerical experiments show that the proposed
scheme yields better approximation than other fifth-order schemes,
especially to the numerical problems which contain discontinuities
while keeping essentially non-oscillatory performance.

\paragraph{Acknowledgements: }

The authors are gratefully thank to the National Board for Higher
Mathematics (NBHM) for the financial support by the Grant Reference
No. $2/48(4)/2014$/NBHM/R \& D II/$/14356$ through the Department
of Atomic Energy (DAE), India.


\begin{thebibliography}{10}
\bibitem[1]{Acker}F. Acker, R. B. de. Borges, and B. Costa. An improved
WENO-Z scheme. J. Comput. Phys., 313:726-753, 2016.

\bibitem[2]{balsaraandshu8} D.S. Balsara, and C.W. Shu. Monotonicity
preserving weighted essentially non-oscillatory schemes with increasingly
high order of accuracy. J. Comput. Phys., 160:405-452, 2000 .

\bibitem[3]{Borgescaramona10} R. Borges, M. Carmona, B. Costa, and
W.S. Don. An improved weighted essentially non-oscillatory scheme
for hyperbolic conservation laws. J. Comput. Phys., 227:3191-3211,
2008.

\bibitem[4]{Castroetal}M. Castro, B. Costa, and W.S. Don. High
order weighted essentially non-oscillatory WENO-Z schemes for hyperbolic
conservation laws. J. Comput. Phys., 230:766-792, 2011.

\bibitem[5]{Pfanetal14}P. Fan, Y. Shen, B. Tian, and C. Yang.
A new smoothness indicator for improving the weighted essentially
non-oscillatory scheme. J. Comput. Phys., 269:329-354, 2014.

\bibitem[6]{Pfan15}P. Fan P. High order weighted essentially non
oscillatory WENO-schemes for hyperbolic conservation laws. J. Comput.
Phys., 269:355-285, 2014.

\bibitem[7]{Gerolymus} G.A. Gerolymos, D. Senechal, and I. Vallet. Very high order WENO schemes. J. Comput. Phys., 228:8481-8524, 2009.

\bibitem[8]{Gottieledandshu21}S. Gottlieb, and C.W. Shu. Total
variation diminishing Runge-Kutta schemes. Math. of Comput., 67:73-85,
1998.

\bibitem[9]{Haetal13} Y. Ha, C.H. Kim, Y.J. Lee, and J. Yoon. An
improved weighted essentially non-oscillatory scheme with a new smoothness
indicator. J. Comput. Phys., 232:68-86, 2013.

\bibitem[10]{Ha83} A. Harten. High resolution schemes for hyperbolic
conservation laws. J. Comput. Phys., 49:357-393, 1983.

\bibitem[11]{Ha84} A. Harten. On a class of high resolution total-variation-stable
finite-difference schemes. SIAM J. Numer. Anal., 21:1, 1984.

\bibitem[12]{HOEC86}A. Harten, S. Osher, B. Engquist, and S.R. Chakravarthy.
Some results on uniformly high-order accurate essentially non-oscillatory
schemes. App. Numer. Math., 2: 347-377, 1986.

\bibitem[13]{HO87}A. Harten, and S. Osher. Uniformly high-order accurate
non oscillatory schemes, I. SIAM J. Numer. Anal., 24: 279-309, 1987.

\bibitem[14]{HO97}A. Harten, B. Engquist, S. Osher, and S.R. Chakravarthy.
Uniformly high order accurate non-oscillatory schemes, III. J. Comput.
Phys., 131:3-47, 1997.

\bibitem[15]{henrickaslampowers5}A.K. Henrick, T.D. Aslam, and
J.M. Powers. Mapped weighted essentially non-oscillatory schemes:
Achieving optimal order near critical points. J. Comput. Phys., 207:
542-567, 2005.

\bibitem[16]{MWENO-Z} F. Hu F, R. Wang, and X. Chen. A modified fifth-order
WENOZ method for hyperbolic conservation laws. J. Comput. and App.
math., 303: 56-68, 2016.

\bibitem[17]{Jiangandshu7} G.S. Jiang, and C.W. Shu. Efficient
implementation of Weighted ENO schemes. J. Comput. Phys., 126:202-228,
1996.

\bibitem[18]{KIMetal}C.H. Kim, Y. Ha, and J. Yoon. Modified Non-linear
Weights for Fifth-Order Weighted Essentially Non-oscillatory Schemes.
J. Sci. Comput., 67:299-323, 2016.

\bibitem[19]{XDLiu}X.D. Liu, S. Osher, and T. Chan. Weighted Essentially
non-oscillatory schemes. J. Comput. Phys., 115:200-212, 1994.

\bibitem[20]{OC84}S. Osher, and S.R. Chakravarthy. High resolution
schemes and the entropy condition. SIAM J. Numer. Anal., 21:5, 1984.

\bibitem[21]{Rinne} C.W. Schulz-Rinne, J.P. Collins, and H.M. Glaz.
Numerical solution of the riemann problem for two-dimensional gas
dynamics. SIAM J. Sci. Compu., 14:1394-1414, 1993 .

\bibitem[22]{shunotes}C.W. Shu. Essentially non-oscillatory and
weighted essentially non-oscillatory schemes for hyperbolic conservation
laws. In Advanced numerical approximation of nonlinear hyperbolic
equations. Lecture Notes in Mathematics, Berlin, Springer-Verlag ,
1697: 325-432, 1993.

\bibitem[23]{cwshu16}C.W. Shu. High order weighted essentially
non oscillatory schemes for convection dominated problems. SIAM Review,
51:82-126, 2009.

\bibitem[24]{osherSHU}C.W. Shu, and S. Osher. Efficient implementation
of essentially non-oscillatory shock-capturing schemes. J. Comput.
Phys., 77: 439-471, 1988.

\bibitem[25]{Shu-osher1}C.W. Shu, and S. Osher. Efficient implementation
of essentially non-oscillatory shock-capturing schemes II. J. Comput.
Phys., 83:32-78, 1989.

\bibitem[26]{Sod}G.A. Sod. A survey of several finite difference
methods for systems of nonlinear hyperbolic conservation laws. J.
Comput. Phys., 107:1-31, 1978.

\bibitem[27]{woodward}P. Woodward, and P. Colella. The numerical
simulation of two-dimensional fluid flow with strong shocks. J. Comput.
Phys., 54:115-173, 1984.
\bibitem[28]{young}Y.N. Young, H. Tufo, A. Dubey,and R. Rosner. On the
miscible Rayleigh-Taylor instability:two and three dimensions. J.
Fluid Mech. 447:377-408, 2001.
\end{thebibliography}
\end{document}